\documentclass[letterpaper,11pt]{amsart}
\usepackage{latexsym,array,delarray,amsthm,amssymb,epsfig,amsmath,blkarray,tikz}
\usepackage{url}
\usepackage{tkz-euclide,subfigure}

\usepackage{longtable}
\usetikzlibrary{decorations.markings}
\usetikzlibrary{arrows}

\usetikzlibrary{arrows.meta}
\usepackage{ytableau}

\theoremstyle{plain}
\newtheorem{thm}{Theorem}[section]

\newtheorem{prop}[thm]{Proposition}

\newtheorem*{thm*}{Theorem}
\newtheorem*{lemma*}{Lemma}
\newtheorem*{prop*}{Proposition}
\newtheorem*{cor*}{Corollary}
\newtheorem*{conj*}{Conjecture}
\newtheorem*{notation}{Notation}

\theoremstyle{definition}
\newtheorem{defn}[thm]{Definition}

\newtheorem{ex}[thm]{Example}

\newtheorem{alg}[thm]{Algorithm}

\theoremstyle{remark}
\newtheorem*{rmk}{Remark}

\newcommand{\zz}{\mathbb{Z}}
\newcommand{\nn}{\mathbb{N}}

\newcommand{\cc}{\mathbb{C}}

\newcommand{\bfa}{\mathbf{a}}

\newcommand{\bfe}{\mathbf{e}}

\newcommand{\bfh}{\mathbf{h}}

\newcommand{\bfu}{\mathbf{u}}

\newcommand{\bfx}{\mathbf{x}}
\newcommand{\bfy}{\mathbf{y}}

\newcommand{\cg}{\mathcal{G}}
\newcommand{\cm}{\mathcal{M}}
\newcommand{\cn}{\mathcal{N}}
\newcommand{\cp}{\mathcal{P}}
\newcommand{\cj}{\mathcal{J}}

\newcommand{\ind}{\mbox{$\perp \kern-5.5pt \perp$}}

\begin{document}
\title[Input-output equations for linear compartmental models]{Input-output equations for linear compartmental models with multiple outputs}
\author{Francisco Ponce Carrion and Seth Sullivant}
\address{Department of Mathematical Sciences \\
Worcester Polytechnic  Institute,  Worcester, MA 01609}
\address{Department of Mathematics \\
North Carolina State University, Raleigh, NC, 27695}
\email{fcarrion@wpi.edu}
\email{smsulli2@ncsu.edu}
\maketitle

\begin{abstract}
Linear compartmental models are linear systems of ODEs that have an associated graphical structure which is often obtained from biological models. Structural identifiability of a model assesses whether one can recover the unknown parameters of a model from input-output data. In order to address this question, we study the set of input-output equations of a linear compartmental model. In this work, we provide an algorithm to find input-output equations where the user can specify what are the orders of the derivatives of the input-output variables that appear. We then study which of the lists of prescribed orders are viable depending on the graphical structure that defines the linear compartmental model.
\end{abstract}

\section{Introduction}
Linear compartmental models (LCMs) are systems of ODEs that have an associated graphical structure, say $G$. This structure gives information about the parameters, the state variables, the input variables, and the output variables. The graphical structure $G$ is often obtained from biological models in pharmacology, ecology, or disease modeling \cite{compartmentPharma,compartmentEcology,compartmentDisease}. In many of these instances, the unknown parameters are physically meaningful, so when doing parameter estimation, it is fundamental to know whether the recovered parameters that yield the desired dynamics are unique, given data, or if there are multiple parameter values that can model the same observations.

Structural identifiability of a model addresses the question of whether the unknown parameters of the model can be inferred from perfect input-output data. One way to tackle this question is to take our system of equations and perform differential elimination to find relationships only involving input-output variables and their derivatives. We then analyze the coefficients, which will be polynomials in the unknown parameters, to determine which of the parameters can be recovered uniquely \cite{treeform, nikkicashIndist, nikkirepar}.

Obtaining input-output equations is a fundamental part of this process
of determining if a model is structurally identifiable. 
Past work has been done describing input-output equations for these models using  subgraphs of graphs associated to $G$ \cite{treeform,nikkirepar} where the
model has a single output. 
When applied to models with multiple outputs, this approach
will produce input-output equations that 
such that each input-output equation is of fixed order, and depends on a single output and a single input variable, so no interactions between different inputs-outputs can be obtained. 
Furthermore, 
those past results are limited to models where $G$ is strongly connected.

In this work, we provide an algorithm to compute input-output equations of prescribed orders for observable systems, 
with an arbitrary sized set of output variables.  In contrast
to the previous work, this can produce a set of input-output equations of lower
order, and which contain multiple outputs in the same equation.
Then we study which prescribed orders of input-output equations are attainable from a given problem.  

Having lower-order input-output equations is useful in practice because it is difficult
to accurately estimate higher order derivatives from data.
Furthermore, as an application of our results we give an example of a 
linear compartmental model for which the graph  is not strongly connected and for which the single-output input-output equations fail to give accurate identifiability results (Example
\ref{ex:nonident}).

The technique for finding the input-output equations for the linear compartmental models has its basis in
differential algebra techniques.  However, the linear structure
of the models allows us to convert problems of differential
elimination into problems of linear algebra on associated matrices of
coefficients.  We relate the calculation of the input-output equations to
the calculation of a reduced row echelon form of an associated matrix.
Whether or not a prescribed order for input-output equations is attainable
is encoded in the invertibility of a certain square submatrix.
Since entries in these matrices are related to paths in an associated 
graph, we approach these problems using combinatorial techniques.

The paper is structured as follows:
In Section \ref{sec:prelim}, we provide the differential algebra tools that we need for the treatment of this problem, following the paper of Dong, Goodbrake, Harrington, and Pogudin \cite{pogudin}. In Section \ref{sec:alg}, we show Algorithm \ref{alg} that takes a linear compartmental model and a valid list of prescribed orders, and returns a set of input-output equations, with the appropriate orders. In Section \ref{sec:correct}, we study the correctness of Algorithm \ref{alg}. In Section \ref{sec:invertibility}, we study the possible valid lists of prescribed orders for input-output equations, by analyzing a determinant whose terms are related to subgraphs in $G$.  We give a necessary condition for the attainability
of a particular set of orders, as an application of Menger's theorem related
to flows in graphs.

\section{Differential Algebra Preliminaries}\label{sec:prelim}

In this section, we begin by formally defining a linear compartmental model, and then we discuss the necessary tools from differential algebra in our treatment of linear compartmental models. Specifically, we go over the concepts of differential rings, differential ideals, and a specific technique of differential elimination. These tools will be useful as they allow us to encode a linear compartmental model in a differential ideal of a differential ring in a similar fashion in which a linear system of polynomials can be encoded in an ideal of a polynomial ring. Furthermore, the language of differential algebra allows us to give a differential-algebraic definition of  identifiability. Finally, we go over differential elimination using parametric profiles and projections, as the problem of finding input-output equations is equivalent to finding generators for a specific elimination ideal.

Linear compartmental models allow us to represent biological systems governed by ODEs by using an underlying directed graph. In this graph, each vertex represents a compartment, each directed edge $i\leftarrow j$ between compartments indicates a transfer from $j$ to $i$ with a transfer rate $a_{ij}$, and if compartment $i$ has a rate of flow that is leaving the system, we say that compartment $i$ has a leak and denote that rate by $a_{0i}$. Each compartment $i$ is identified with the state variables $x_i$. The $x_i$ that can be measured are called the output variables, and are denoted with an undirected edge on that compartment. We may also have inputs in our system, which will be denoted by an incoming edge pointing toward a compartment. We encode all of this information in a tuple $(G,In,Out,Leak)$ as follows:
\begin{defn}
	A \emph{linear compartmental model (LCM)} is a tuple $\mathcal{M} = (G, In, Out, Leak)$ where $G = (V_G,E_G)$ is a directed graph, and $In,Out,Leak \subseteq V_G$. We will denote $|V_G| = n$. Every linear compartmental model has the following associated dynamical system:
	$$\Sigma_\mathcal{M} = \begin{cases}
	\bfx'(t) &= A\bfx(t)+\bfu(t)\\
	y_i(t) &= x_i(t)\quad \mbox{ for } i\in Out\\
	y_i(t) &= 0 \quad else\\
	u_i(t) &= 0 \quad \mbox{ for } i\notin In
	\end{cases}$$
	where
	\begin{enumerate}
		\item $A_{i,j} = \begin{cases}  
-a_{0i} - \sum_{k\colon k \leftarrow i\in E_G} a_{ki}  & i=j,\  i\in Leak, \\
-\sum_{k\colon k \leftarrow i  \in E_G} a_{ki} & i=j, \ i \not\in Leak, \\
a_{ij} 		& i\neq j,~  (j,i) \in E_G, \\
0 			&  i\neq j,~  (j,i) \notin E_G. 
\end{cases}$
		\item $\bfx(t) = \left(x_1(t),\ldots, x_n(t)\right)$ is the vector of state variables.
		\item $\bfy(t) = \left(y_1(t),\ldots, y_n(t)\right)$ is the vector of output variables.
		\item $\bfu(t) = \left(u_1(t),\ldots, u_n(t)\right)$ is the vector of input variables.
	\end{enumerate}	
\end{defn}

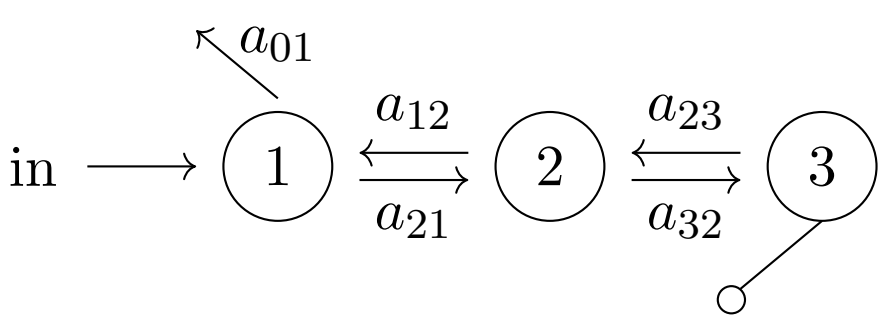
\begin{figure}
	\begin{center}
	\begin{tikzpicture}[scale=1.8]
 	\draw (0,0) circle (0.2);	
 	\draw (1,0) circle (0.2);	
 	\draw (2,0) circle (0.2);	
    	\node[] at (0, 0) {1};
    	\node[] at (1, 0) {$2$};
    	\node[] at (2, 0) {$3$};
	 \draw[<-] (0.3, 0.05) -- (0.7, 0.05);	
	 \draw[->] (0.3, -0.05) -- (0.7, -0.05);	
	 \draw[<-] (1.3, 0.05) -- (1.7, 0.05);	
	 \draw[->] (1.3, -0.05) -- (1.7, -0.05);	
   	 \node[] at (0.5, 0.2) {$a_{12}$};
   	 \node[] at (1.5, 0.2) {$a_{23}$};
   	 \node[] at (0.5, -0.2) {$a_{21}$};
   	 \node[] at (1.5, -0.2) {$a_{32}$};
 	\draw (5/3,-.49) circle (0.05);	
	 \draw[-] (2, -.2 ) -- (1.7, -.45);	
	 \draw[->] (-0.7, 0) -- (-0.3, 0);	
   	 \node[] at (-.9, 0) {in};
	 \draw[->] (0, .25) -- (-0.3, .5);	
   	 \node[] at (0, 0.45) {$a_{01}$};
	\end{tikzpicture}
\end{center}

\caption{\label{fig:catenary} A catenary linear compartment model with one input and one output.}
\end{figure}

\begin{ex}
	Consider the linear compartmental model $\cm = (G, In,Out,Leak)$ given by
	the graph in Figure  \ref{fig:catenary},.

The system of ODEs in this case will be
\begin{equation*}
	\begin{aligned}
		\begin{pmatrix}
			x_1'\\x_2'\\x_3'
		\end{pmatrix} &= \begin{pmatrix}
			-a_{01} - a_{21} & a_{12} & 0 \\
			a_{21} & -a_{12}-a_{32} & a_{23}\\
			0 & a_{32} & - a_{23}
		\end{pmatrix}\begin{pmatrix}
			x_1\\x_2\\x_3
		\end{pmatrix} + \begin{pmatrix}
			u_1 \\ 0\\0
		\end{pmatrix},\\
		\begin{pmatrix}
			y_1\\y_2\\y_3
		\end{pmatrix} &= \begin{pmatrix}
			0\\0\\x_3
		\end{pmatrix}.
	\end{aligned}
\end{equation*}
\end{ex}

We are interested in studying the systems of differential equations of linear compartmental models from an algebraic perspective. We study such systems by using tools from \emph{differential algebra}. The main objects we will use from this area are: 
\begin{enumerate}
	\item Differential rings: these will generally be polynomial rings in variables $\bfx$, $\bfy$, $\bfu$ (and their derivatives) with coefficients in $\cc (a_{ij}: 1\leq i,j  \leq n)$. The derivation is classical differentiation with respect to a time variable, $t$, where we treat $a_{ij}$ as constants. These rings are where our systems of ODEs live. 
	\item Differential ideals: we will encode the system of ODEs in an ideal which is closed under differentiation.
	\item Projection techniques: input-output equations are the result of elimination. We show a specific method using projections to eliminate variables that guarantees our resulting equations are sufficient to conclude about identifiability.
\end{enumerate}

\begin{notation}\
\begin{enumerate}
	\item Let $x$ be an element of a differential ring and $h\in \zz_{\geq0}$. 
	\begin{align*}
		x^{(<h)} &= \left(x,x', \ldots, x^{(h-1)}\right),\\
		x^{(\infty)} &= {\left(x,x',\ldots\right)}.
	\end{align*}
	\item Let $\bfx = (x_1,\ldots,x_n)$ with $x_i$ in a differential ring, $h\in \zz_{\geq0} \cup \infty$, and $\bfh = (h_1,\ldots, h_n)$ with $h_i \in \zz_{\geq0} \cup \infty$, then 
\begin{align*}
	\bfx^{(< h)} &= (x_1^{(<h)},\ldots, x_n^{(<h)}),\\
	\bfx^{(< \bfh)} &= (x_1^{(<h_1)},\ldots, x_n^{(<h_n)}).
\end{align*}
\item Let  $[n] = \{1,\ldots, n\}$, and $H_1,H_2 \subseteq [n]$. Then we define the submatrix $A_{H_1,H_2}$ of $A$ as follows:
\begin{align*}
	A_{H_1,H_2} = (a_{ij}:\; i \in H_1,\; j \in H_2) 
\end{align*}
If $H_1 = \{k\}$ for some $k\in [n]$, then we use the shorthand $A_{k,H_2}$, similarly for $H_2$.
\item Let $\cc(\bfa) = \cc (a_{ij}: 1\leq i,j  \leq n)$.
\end{enumerate}
\end{notation}

Now we are ready to state the differential ideal that encodes the system of ODEs associated to a linear compartmental model. 
\begin{defn}
	Let $\cm=(G, In, Out, Leak)$ be a linear compartmental model. We define the \emph{differential ideal associated to $\cm$, $\cj_\cm$}, as 
	$$\cj_{\mathcal{M}} = \langle \left(x_i'-A_{i,[n]}\bfx\right)^{(\infty)}, \left(y_j-x_j\right)^{(\infty)}: 1\leq i \leq n, j \in Out \rangle\subseteq \cc(\bfa)[x^{(\infty)},y^{(\infty)},u^{(\infty)}].$$
\end{defn}
Every element of $\cj_{\cm}$ vanishes on every analytic solution of its associated ODE system, $\Sigma_\mathcal{M}$. Furthermore, this ideal is prime (see \cite{sian}), and it is the smallest differential ideal containing our defining equations (see \cite{pogudin}). 

We are interested in studying input-output equations due to their usefulness in deciding whether a model is identifiable or not. Next, we give algebraic definitions of the different types of identifiability.

\begin{defn}\
\begin{enumerate}
	\item A rational function $f \in \cc(\bfa)$ is \emph{globally single experiment (SE) identifiable} if there exist $P_0,P_1\in \cc[\bfy^{(\infty)},\bfu^{(\infty)}]$ such that $P_1\notin \cj_{\cm}$ and $P_1 f-P_0 \in \cj_{\cm}$.
	
	\item $f$ is \emph{locally SE identifiable} if there exist $s\in \nn$ and $P_0,\ldots, P_s \in \cc[y^{(\infty)},u^{(\infty)}]$ such that 
	$$P_s f^s + P_{s-1} f^{s-1} +\ldots + P_0 \in \cj_{\cm}.$$
	\item Let $r\in \nn$, the \emph{$r-$fold replica} of the model associated to $\cm$ is given by
	$$\Sigma_\mathcal{M}^r = \begin{cases}
	\bfx_k' = A\bfx_k+\bfu_k & k \in [r]\\
	y_{ik} = x_{ik}\quad \mbox{ for } i\in Out &k \in [r]\\
	y_{ik} = 0 \quad \quad \mbox{ for } i\notin Out &k\in [r]\\
	u_{ik} = 0 \quad \mbox{ for } i\notin In & k\in [r]
	\end{cases}$$
	where $\bfx_k = (x_{1k},\ldots, x_{nk})$, and similarly for $\bfy_k, \bfu_k$.  
	\item A rational function $f \in \cc(\bfa)$ is \emph{globally (resp. locally) multi-experiment identifiable} if there exists $r\in\nn$ such that $f$ is globally (resp. locally) SE identifiable in $\Sigma_\cm^r$.
\end{enumerate}	
\end{defn}

\begin{defn}
	Let $J\subseteq \cc(\bfa)[\bfx^{(\infty)}]$ be a differential ideal over $\cc(\bfa)$. The \emph{field of definition} of $J$ is the smallest subfield $L\subseteq \cc(\bfa)$ such that $J$ is generated by $J\cap L[\bfx^{(\infty)}]$.
\end{defn}

\begin{thm}[\cite{pogudin}]\label{thm:gleb}
	The field of multi-experimental identifiable functions is generated over $\cc$ by the field of definition of 
	$$\cj_{\cm} \cap \cc(\bfa)[\bfy^{(\infty)},\bfu^{(\infty)}].$$
\end{thm}

Theorem \ref{thm:gleb} allows us to answer a more general question than whether a model is ME identifiable or not. If we can compute generators of the field of definition of $\cj_{\cm} \cap \cc(\bfa)[\bfy^{(\infty)},\bfu^{(\infty)}]$, we can then decide which parameters (or functions of them) are identifiable. Indeed, if we compute the generators of the elimination ideal $\cj_{\cm} \cap \cc(\bfa)[\bfy^{(\infty)},\bfu^{(\infty)}]$ in a specific way, there is an algorithm to compute the generators for the field of definition of the elimination ideal (see \cite{pogudin}, Algorithm 4.5). Specifically, we need the generators for the elimination ideal (input-output equations) to be computed through the use of parametric profiles and projections, which we will discuss next.

\begin{defn}
	For a prime differential ideal $P$, we say $\bfh = (h_1,\ldots, h_n)$ is a \emph{parametric profile} for $P$ if the images of $\bfx^{(<\bfh)}$ form a transcendence basis of $\cc(\bfa)[\bfx^{(\infty)}]/P$.
\end{defn}

\begin{defn}
	For a prime differential ideal $P \subset  K[\bfx ^{(\infty )}]$ and
its parametric profile $\bfh$, a tuple of differential polynomials $(f_1,\ldots, f_n)$ is called the \emph{projection
corresponding to $\bfh$}  if, for every $1\leq i\leq n$, the following hold:
\begin{enumerate}
	\item If $h_i =\infty$, then $f_i =0$.
	\item If $h_i <\infty$, then $f_i$ is the generator of $P \cap \cc(\bfa)[\bfx^{(<\bfh +\bfe_i)}]$.
\end{enumerate} 
\end{defn}
\begin{rmk}
	For a linear compartmental model, $\cm$, the ideal $\cj_{\cm}$ admits the following standard parametric profile:
	$$\bfh_s = (\underbrace{1,\ldots, 1}_{n\mbox{ times}},\underbrace{0, \ldots,0}_{\left|Out\right|\mbox{ times}}, \underbrace{\infty, \ldots, \infty}_{\left|In\right|\mbox{ times}})$$
with corresponding projections 
	$$(x_1'-A_{1,[n]}\bfx, \ldots, x_n'-A_{n,[n]}\bfx,y_{j_1}-x_{j_1}, \ldots, y_{j_{k}}-x_{j_{k}} , \underbrace{0, \ldots, 0}_{\left|In\right|\mbox{ times}})$$
where $Out = \{j_1,\ldots, j_k\}$.
\end{rmk}

An ideal admits many different representations by projections with different profiles. The Socoban algorithm described in \cite{pogudin} describes a way to "move" between representations to get a desired one. The key idea is that we want to move the 1s in the profile corresponding to the state variables towards the slots corresponding to outputs and inputs (these projections will then be our input-output equations). We are interested in parametric profiles, $\bfh_{\lambda}$, where we move all the 1s from the state variables into the output variables, namely:
\begin{equation}\label{eq:tprof}
	\bfh_{\lambda} = (\underbrace{0,\ldots, 0}_{n\mbox{ times}},\underbrace{\lambda_1, \ldots,\lambda_{\left|Out\right|}}_{\left|Out\right|\mbox{ times}}, \underbrace{\infty, \ldots, \infty}_{\left|In\right|\mbox{ times}}), \mbox{ where } \lambda_1 + \cdots + \lambda_{\left|Out\right|}=n.
\end{equation} 
The corresponding projections to the output variables will be input-output equations of orders $\lambda_i$, for $i \in Out$. The zeros in the state variables imply that the system $\Sigma_\cm$ is \textit{observable}. This is because the projection corresponding to the $i$-th state variables will only contain one $x_i$, and input-output variables, for each $i \in [n]$. 

Work done by Dong, Goodbrake, Harrington, and Pogudin \cite{pogudin}, shows how to compute generators for the field of ME-identifiable functions given a parametric profile and its projections. This leads to the main question we address in this work. Given a linear compartmental model, $\mathcal{M} = (G, In, Out, Leak)$, what are the distinct parametric profiles and projections one can obtain?

\section{Finding input-output equations using row elimination}\label{sec:alg}

In this section, we provide an algorithm to compute input-output equations of a linear compartmental model that takes as input a linear compartmental model and a parametric profile $\bfh_\lambda$. Providing $\bfh_\lambda$ as an input to our algorithm allows us to decide the order of the derivatives of the input-output variables that appear in our input-output equations. In addition, this elimination algorithm takes advantage of the linear structure of the equations that define a linear compartmental model as we perform differential elimination by simply using row elimination on an appropriate matrix.

We are interested in studying parametric profiles of the form of \eqref{eq:tprof}. To simplify notation, we will make the following assumptions throughout this work:
	\begin{enumerate}
		\item Every compartment of $\mathcal{M}$ has a leak, so we can treat the diagonal entries of $A$ as independent coefficients $a_{ii}$, $i\in [n]$.
		\item $Out = \{n-k+1,\ldots, n\}$ for some $1\leq k<n$.
		\item $\lambda = (\underbrace{0,\ldots, 0}_{n-k\mbox{ times}},\lambda_{n-k+1},\ldots, \lambda_n)$
		\item $0< \lambda_{n-k}\leq \ldots\leq \lambda_{n}$.
		\item $\bfh_{\lambda} = (\underbrace{0,\ldots, 0}_{n\mbox{ times}},\lambda_{n-k+1}, \ldots,\lambda_{n}, \underbrace{\infty, \ldots, \infty}_{\left|In\right|\mbox{ times}})$.
	\end{enumerate} 
It follows from the observation that if we have a leak in each compartment, we can use the change of variables $\left(-a_{0i} - \sum_{k\colon k \leftarrow i\in E_G} a_{ki}\right) \rightarrow a_{ii}$ for each $i$. Assumptions 2-5 are mainly so that we can identify a target parametric profile $\bfh_\lambda$ to a non-decreasing list $\lambda$, whose zero entries are state variables, and non-zero entries are output variables. 

Let us consider a linear compartmental model $\cm=(G, In, Out, Leak)$ with $n$ state variables. We will be looking at the system defined by
\begin{equation*}
	f= \bfx'-A\bfx-\bfu.
\end{equation*}
\begin{notation}\
	\begin{enumerate}
		\item $\gamma = \max(\lambda)$.
		\item $D = diagonal\left(\sum_{i \in In} e_i\right)$.
		\item $S^{\bullet}_{l}= \left\{ i \in [n]: \lambda_i \bullet l \right\}$, where $\bullet \in \{<,>,\leq, \geq\}$.
		\item $S_{l} = \left\{ i \in [n]: \lambda_i = l \right\} $
		\item $I_m$ denotes the $m\times m$ identity matrix, for any $m\in \nn$.
	\end{enumerate}
\end{notation}

%

\begin{alg}\label{alg}\

\noindent In: Linear compartmental model $\cm=(G, In, Out, Leak)$, and parametric profile $\bfh_{\lambda}$.\\
Out: Input-output equations whose variables have orders specified in $\bfh_\lambda$.
\begin{enumerate}
	\item Prolongate the system $f$ into $F$, where  $$F = \left( f_i^{(\gamma)}: \lambda_i = \gamma\right) \cup \left(f^{(<\gamma)}\right). $$
	\item Write down the system $F$ in matrix form:
	\begin{equation*}
		M = \bordermatrix{
	&\bfx_{S_\gamma}^{(\gamma)}&\bfx^{(\gamma-1)} & \bfx^{(\gamma-2)} \cdots & \bfx^{(1)}& \bfx^{(0)} & \bfu^{(\gamma-1)} & \bfu^{(\gamma-2)} & \cdots & \bfu^{(0)}\cr
	& -I_{\left|S_\gamma\right|} & A_{S_\gamma,[n]}& & & & D_{S_\gamma,[n]} & & &\cr
	&& -I_n & A& & & & D &  &\cr
	&& & & -I_n& A & &  &  &D\cr
	}
	\end{equation*}
	\item Take a monomial ordering satisfying $\bfu^{(0)}<\cdots < \bfu^{(\gamma-1)} <\bfx_{S^{>}_{0}}^{(0)} < \bfx_{S^{>}_{ 1}}^{(1)}<\cdots < \bfx_{S^{>}_{ \gamma -1}}^{(\gamma-1)} < \bfx_{S_{0}}^{(0)} < \bfx_{S^{\leq}_{ 1}}^{(1)}<\cdots < \bfx_{S^{\leq}_{ \gamma -1}}^{(\gamma-1)}< \bfx_{S_{\gamma}}^{(\gamma)}$. Then define $\bar{M}$ as the reordering of the columns of $M$ in descending order. 
	\item We obtain the following block structure:
	\begin{equation*}
		\bar{M} = \bordermatrix{
	&\bfx_{S_\gamma}^{(\gamma)}&\bfx_{S^{\leq}_{ \gamma-1}}^{(\gamma-1)} & \bfx_{S^{\leq}_{ \gamma -2}}^{(\gamma-2)}& \cdots & \bfx_{S^{\leq}_{ 1}}^{(1)}& \bfx_{S_{0}}^{(0)} & \bfx_{S^{>}_{ \gamma -1}}^{(\gamma-1)}& \cdots & \bfx_{S^{>}_{0}}^{(0)} & \bfu^{(\leq\gamma-1)} \cr
	& -I_{\left|S_\gamma\right|} & A_{S_{\gamma}, S^{\leq}_{ \gamma-1}}& & && & & & &\cr
	&& -I_{\left|S^{\leq}_{ \gamma-1}\right|} & A_{[n],S^{\leq}_{ \gamma -2}}&\ddots & & & & \star &  &\star\star\cr
	&& && & -I_{\left|S^{\leq}_{ 1}\right|}& A_{[n],S_{0}} & &  &  &\cr
	}
	\end{equation*}
	\item Compute the reduced row echelon form (RREF) of $\bar{M}$. Input-output equations are the pivots corresponding to $\bfx_{Out}^{\lambda}$.
\end{enumerate}	
\end{alg}
\begin{rmk}
Note that we are assuming that the input $\bfh_\lambda$ is a parametric profile for the ideal $\cj_\cm$. In Section \ref{sec:correct}, we show how this plays a fundamental role in proving the correctness of the output of Algorithm \ref{alg}. 
\end{rmk}

\begin{defn}
	Let $\cm=(G, In, Out, Leak)$ be a linear compartmental model on $n$ state variables, $\bfh_\lambda$ be a parametric profile and let $\bar{M}$ be defined as in step (4) of Algorithm \ref{alg}. We define \emph{$M_\lambda$} to be the first $(\gamma-1)n$ columns of $\bar{M}$ (the columns of $\bar{M}$ before the stars), i.e.:
	$$M_\lambda = \bordermatrix{
	&\bfx_{S_\gamma}^{(\gamma)}&\bfx_{S^{\leq}_{ \gamma-1}}^{(\gamma-1)} & \bfx_{S^{\leq}_{ \gamma -2}}^{(\gamma-2)}& \cdots & \bfx_{S^{\leq}_{ 1}}^{(1)}& \bfx_{S_{0}}^{(0)} \cr
	& -I_{\left|S_\gamma\right|} & A_{S_{\gamma}, S^{\leq}_{ \gamma-1}}& & && \cr
	&& -I_{\left|S^{\leq}_{ \gamma-1}\right|} & A_{[n],S^{\leq}_{ \gamma -2}}&\ddots & &\cr
	&& && & -I_{\left|S^{\leq}_{ 1}\right|}& A_{[n],S_{0}}\cr
	}. $$
\end{defn}
The question of the existence of input-output equations with a given parametric profile $\bfh_\lambda$ is then related to the reduced-row echelon form of $\bar{M}$, which we study by analyzing the invertibility of $M_\lambda$. First, let us see an example where we show how Algorithm \ref{alg} works.
\begin{ex}\label{ex:ex1} Let us consider the linear compartmental model $\mathcal{M} =  (G, In, Out, Leak)$ with $In = \{1\}$, $Out = \{2,3\}$ and $Leak = [3]$, and $G$ as illustrated in Figure \ref{fig:ex1}.
\begin{figure}
\begin{center}
	\begin{tikzpicture}[scale=1.8]
 	\draw (0,0) circle (0.2);	
 	\draw (1,0) circle (0.2);	
 	\draw (2,0) circle (0.2);	
    	\node[] at (0, 0) {1};
    	\node[] at (1, 0) {$3$};
    	\node[] at (2, 0) {$2$};
	 \draw[<-] (0.3, 0.05) -- (0.7, 0.05);	
	 \draw[->] (0.3, -0.05) -- (0.7, -0.05);	
	 \draw[<-] (1.3, 0.05) -- (1.7, 0.05);	
	 \draw[->] (1.3, -0.05) -- (1.7, -0.05);	
   	 \node[] at (0.5, 0.2) {$a_{13}$};
   	 \node[] at (1.5, 0.2) {$a_{32}$};
   	 \node[] at (0.5, -0.2) {$a_{31}$};
   	 \node[] at (1.5, -0.2) {$a_{23}$};
 	\draw (5/3,-.49) circle (0.05);	
	 \draw[-] (2, -.2 ) -- (1.7, -.45);	
 	\draw (2/3,-.49) circle (0.05);	
	 \draw[-] (1, -.2 ) -- (0.7, -.45);	
	 \draw[->] (-0.7, 0) -- (-0.3, 0);	
   	 \node[] at (-.9, 0) {in};
	\end{tikzpicture}
\end{center}
\caption{Linear compartmental model in Example \ref{ex:ex1}}\label{fig:ex1}
\end{figure}
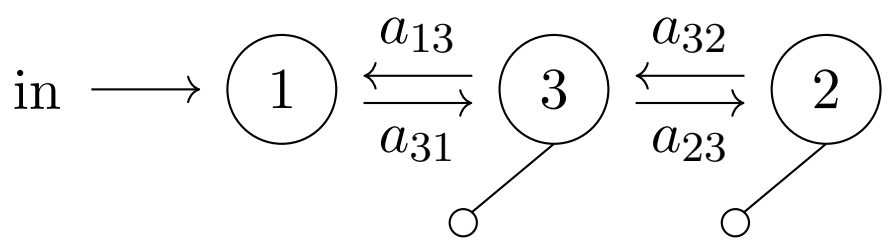

This corresponds to the system \begin{equation*}
	\begin{aligned}
		x_1' &= a_{11}x_1 + a_{13}x_3 + u_1\\
		x_2' &= a_{22}x_2 + a_{23}x_3\\
		x_3' &= a_{31}x_1 + a_{32}x_2 + a_{33}x_3\\
		y_2&=x_2\\
		y_3&=x_3
	\end{aligned}
\end{equation*}

First, we can rewrite our original system as

\begin{equation*}
	\bordermatrix{
	&x_1'&x_2'&x_3'&x_1&x_2&x_3 &u_1\cr
	&-1 &0&0 & a_{11} & 0 & a_{13} &1 \cr
	&0 &-1&0 & 0 & a_{22} & a_{23} &0 \cr
	&0 &0&-1 & a_{31} & a_{32} & a_{33}&0 \cr
	}
\end{equation*}

We want to get to the profile associated to  $\lambda = (0,1,2)$. We will augment our matrix to have $x_3''$ as follows:

\begin{equation*}
	M=\bordermatrix{
	&x_3''& x_1'&x_2'&x_3'&x_1&x_2&x_3&u_1\cr
	&-1&a_{31}&a_{32}&a_{33} &0 &0 & 0 &0\cr
	&0&-1 &0&0 & a_{11} & 0 & a_{13} &1\cr
	&0&0 &-1&0 & 0 & a_{22} & a_{23} &0 \cr
	&0&0 &0&-1 & a_{31} & a_{32} & a_{33}&0 \cr
	}
\end{equation*}
We then replace labelings with output variables to see clearly which of the columns correspond to output variables.
\begin{equation*}
	\bordermatrix{
	&y_3''& x_1'&y_2'&y_3'&x_1&y_2&y_3 &u_1\cr
	&-1&a_{31}&a_{32}&a_{33} &0 &0 & 0 &0 \cr
	&0&-1 &0&0 & a_{11} & 0 & a_{13} &1\cr
	&0&0 &-1&0 & 0 & a_{22} & a_{23} &0\cr
	&0&0 &0&-1 & a_{31} & a_{32} & a_{33}&0 \cr
	}
\end{equation*}
Then, we reorganize the columns to obtain $\bar{M}$. The pivots in blue are the variables we want to solve for, and the pivots in red are the variables we want to eliminate.
\begin{equation*}
	\bar{M} = \bordermatrix{
	&\color{blue}y_3''&\color{red}x_1'&\color{blue}y_2'&\color{red}x_1&y_3'&y_2&y_3 &u_1\cr
	&-1&a_{31}&a_{32} & 0      & a_{33} & 0      & 0     &0 \cr
	&0 &-1    &0      & a_{11} & 0      & 0      & a_{13}&	1 \cr
	&0 &0     &-1     & 0      & 0      & a_{22} & a_{23}&0 \cr
	&0 &0     &0      & a_{31} & -1     & a_{32} & a_{33}&0 \cr
	}
\end{equation*}
Note that $$M_\lambda = \bordermatrix{
	&\color{blue}y_3''&\color{red}x_1'&\color{blue}y_2'&\color{red}x_1\cr
	&-1&a_{31}&a_{32} & 0       \cr
	&0 &-1    &0      & a_{11}  \cr
	&0 &0     &-1     & 0       \cr
	&0 &0     &0      & a_{31}  \cr
	}.$$
We then compute the reduced row echelon form (RREF) of $\bar{M}$.
\begin{equation*}
	RREF(\bar{M}) = \bordermatrix{
	&\color{blue}y_3''&\color{red}x_1'&\color{blue}y_2'&\color{red}x_1&y_3'&y_2&y_3 &u_1\cr
&1 & 0 & 0 & 0 & -a_{11}-a_{33} & a_{11}a_{32}-a_{22}a_{32} & a_{11} a_{33}-a_{13} a_{31}-a_{23} a_{32}&-a_{31}\cr
& 0 & 1 & 0 & 0 & -\frac{a_{11}}{a_{31}} & \frac{a_{11} a_{32}}{a_{31}} & \frac{a_{11} a_{33}-a_{13} a_{31}}{a_{31}} &-1\cr
 &0 & 0 & 1 & 0 & 0 & -a_{22} & -a_{23}&0\cr
 &0 & 0 & 0 & 1 & -\frac{1}{a_{31}} & \frac{a_{32}}{a_{31}} & \frac{a_{33}}{a_{31}}&0 \cr
}.
\end{equation*}
%
The input-output equations are then
$$\left(\begin{array}{c}
y_{3}''+ \left(-a_{11}-a_{33}\right) y_{3}'+  \left(a_{11}a_{32}-a_{22}a_{32}\right) y_{2}+\left(a_{11} a_{33}-a_{13} a_{31}-a_{23} a_{32}\right) y_{3}-a_{31}u_1
\\
 y_2'-a_{22} y_{2}-a_{23} y_{3}
\end{array}\right) = \begin{pmatrix}
	0\\0
\end{pmatrix}$$
Note that the last equation, $$x_1 -\frac{1}{a_{31}}y_3' + \frac{a_{32}}{a_{31}} y_2+ \frac{a_{33}}{a_{31}}y_3 =0$$ allows us to express $x_1$ in terms of the output variables, hence showing the model is observable.
\end{ex}
When our model has multiple outputs, and in particular, when the underlying graphs are not strongly connected, the approach of computing each input-output equation separately using Cramer's rule (see \cite{treeform}) can yield inconclusive results for identifiability. Consider the following example.
\begin{ex}\label{ex:nonident}
Let $\mathcal{M} =  (G, In, Out, Leak)$ be the linear compartmental model with $In = \{3\}$, $Out = \{2,3\}$, $Leak = [3]$ and $G$ as illustrated in Figure \ref{fig:nonident}.

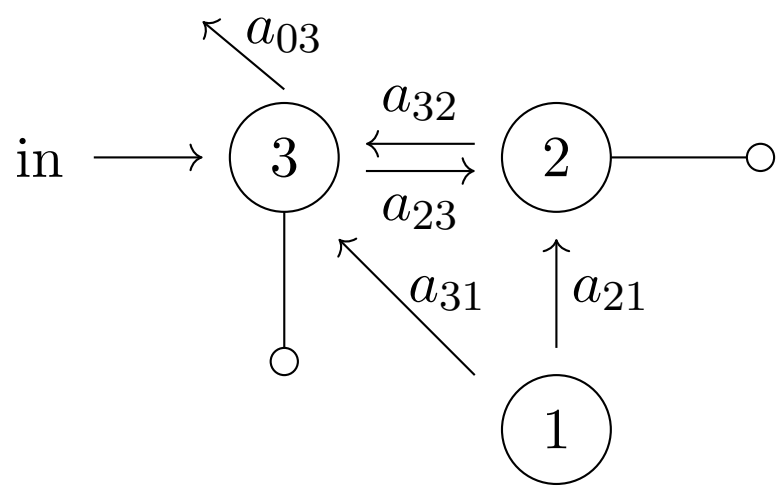
\begin{figure}
\begin{center}
		\begin{tikzpicture}[scale=1.8]
 	\draw (0,0) circle (0.2);	
 	\draw (1,0) circle (0.2);	
 	\draw (1,-1) circle (0.2);	
    	\node[] at (0, 0) {3};
    	\node[] at (1, 0) {$2$};
    	\node[] at (1, -1) {$1$};
	 \draw[<-] (0.3, 0.05) -- (0.7, 0.05);	
	 \draw[->] (0.3, -0.05) -- (0.7, -0.05);	
	 \draw[<-] (0.2, -0.3) -- (0.7, -0.8);
	 \draw[<-] (1, -0.3) -- (1, -0.7);
   	 \node[] at (0.5, 0.2) {$a_{32}$};
   	 \node[] at (0.5, -0.2) {$a_{23}$};
   	 \node[] at (0.6, -0.5) {$a_{31}$};
   	 \node[] at (1.2, -0.5) {$a_{21}$};
 	\draw (1.75,0) circle (0.05);	
	 \draw[-] (1.2, 0 ) -- (1.7, 0);	
	 \draw (0,-0.75) circle (0.05);	
	 \draw[-] (0, -0.2 ) -- (0, -0.7);	
	 \draw[->] (-0.7, 0) -- (-0.3, 0);	
   	 \node[] at (-.9, 0) {in};
	 \draw[->] (0, .25) -- (-0.3, .5);	
   	 \node[] at (0, 0.45) {$a_{03}$};
	\end{tikzpicture}
	\end{center}
	\caption{Linear compartmental model in Example \ref{ex:nonident}}\label{fig:nonident}
	\end{figure}
	Using techniques from \cite{treeform}, we can find the following input-output equations combinatorially:
	\begin{equation*}
		\begin{aligned}
			x_3'''+&\left(\left(-a_{21}-a_{31}\right) a_{32} \left(-a_{03}-a_{23}\right)+\left(-a_{21}-a_{31}\right) a_{23} a_{32}\right) x_{3}\\+ &\left(\left(-a_{32}-a_{03}-a_{23}\right) \left(-a_{21}-a_{31}\right)-a_{32} \left(-a_{03}-a_{23}\right)-a_{23} a_{32}\right) x_3' \\ +&\left(a_{21}+a_{31}+a_{32}+a_{03}+a_{23}\right) x_3''+\left(-a_{21}-a_{31}\right) a_{32} u_{3}+\left(-a_{21}-a_{31}-a_{32}\right) u_3'-u_3''=0,\\
x_2'''+&\left(\left(-a_{21}-a_{31}\right) a_{32} \left(-a_{03}-a_{23}\right)+\left(-a_{21}-a_{31}\right) a_{23} a_{32}\right) x_{2}\\+&\left(\left(-a_{32}-a_{03}-a_{23}\right) \left(-a_{21}-a_{31}\right)-a_{32} \left(-a_{03}-a_{23}\right)-a_{23} a_{32}\right) x_2'\\ +&\left(a_{21}+a_{31}+a_{32}+a_{03}+a_{23}\right) x_2''+\left(-a_{21}-a_{31}\right) a_{23} u_{3}-a_{23}u_3'=0.
		\end{aligned}
	\end{equation*}
	Using the coefficient map method (see \cite{treeform}), and computing the rank of its Jacobian, we would not be able to conclude this model is identifiable through this set of equations. Using our method, we find the following input-output equations corresponding to the profile $\bfh_\lambda$ where $\lambda = (0,1,2)$:
	\begin{equation*}
		\begin{aligned}
		x_{3}^{\prime \prime}+&\frac{\left(a_{32} \left(-a_{21}-a_{31}+a_{32}\right) a_{31}+a_{21} a_{32}^{2}\right)}{a_{31}}x_{2}\\+&\frac{\left(\left(-a_{21}-a_{31}\right) a_{31} \left(-a_{03}-a_{23}\right)+a_{21} a_{32} \left(-a_{03}-a_{23}\right)-a_{23} a_{31} a_{32}\right) }{a_{31}}x_{3}\\+&\frac{\left(\left(a_{21}+a_{31}+a_{03}+a_{23}\right) a_{31}-a_{32} a_{21}\right) }{a_{31}}x_{3}^{\prime}+\frac{\left(\left(-a_{21}-a_{31}\right) a_{31}+a_{32} a_{21}\right)}{a_{31}}u_{3}-u_{3}^{\prime} =0,\\
		x_{2}^{\prime}+&\frac{\left(a_{32} a_{21}+a_{32} a_{31}\right)}{a_{31}}x_{2}+\frac{\left(a_{21} \left(-a_{03}-a_{23}\right)-a_{31} a_{23}\right)}{a_{31}}x_3 -\frac{a_{21} }{a_{31}}x_{3}^{\prime}+\frac{a_{21}}{a_{31}}u_{3}=0.
 		\end{aligned}
	\end{equation*}
	Either using the rank of the Jacobian of the coefficient map, or by computing the generators of the field of definition of the elimination ideal, one can show this model is identifiable using the latter set of equations.
\end{ex}

\section{Correctness of Algorithm \ref{alg}}\label{sec:correct}
Given a linear compartmental model $\cm = (G, In, Out, Leak)$, and a target parametric profile $\bfh_\lambda$, we can construct the matrix $M_\lambda$. In this section, we study $M_\lambda$, and in particular, show how its invertibility is a sufficient condition for $\bfh_\lambda$ to be a parametric profile. We prove this result to show the correctness of Algorithm \ref{alg} in Theorem \ref{thm:correct}.

\begin{prop}
	 Let $\cm$ be a linear compartmental model, and $\bfh_\lambda$ a target parametric profile. Then $M_\lambda$ is a square matrix if and only if $\sum_{i =1}^{n} \lambda_i = n$.
\begin{proof}
Let $n_r$ and $n_c$ be the number of rows and the number of columns of $M_\lambda$, respectively. Note that 
\begin{equation*}
	\begin{aligned}
		n_r &= |S_\gamma| + (\gamma-1)n\\
		n_c &= |S_\gamma| + \sum_{i = 0}^\gamma i |S_{\gamma-i}|.\\
	\end{aligned}
\end{equation*}
Thus, $M_\lambda$ is square if and only if $$\gamma n = n + \sum_{i = 0}^\gamma i |S_{\gamma-i}|.$$
Note that $\sum_{i \in \lambda} i = n$ if and only if $\sum_{i = 0}^\gamma i|S_i|=n$. We look at the
lefthand side and obtain
\begin{align*}
	\gamma n &= \gamma \sum_{i = 0}^\gamma |S_i|\\
	&=  \sum_{i = 0}^\gamma \gamma|S_i| \\
	&= \sum_{i = 0}^\gamma (\gamma-i)|S_i| + \sum_{i = 0}^\gamma i|S_i| \\
	&= \sum_{i = 0}^\gamma i|S_{\gamma - i}| + n.   \qedhere
\end{align*}
\end{proof}
\end{prop}

\begin{ex}
We can picture $\lambda$ as a tableau, so that the equality $\gamma n = n + \sum_{i = 0}^\gamma i |S_{\gamma-i}|$ can be seen as follows. Let $\lambda = (0,0,0,1,1,2,3)$, then we can associate the tableau pictured in Figure \ref{fig:tableau}. Note that in green, we have the tableau corresponding to the partition $(1,1,2,3)$, from which we can construct the corresponding complementary tableau (pictured in white). The complementary tableau to $\lambda = (0,0,0,1,1,2,3)$ would be $\lambda' = (3,3,3,2,2,1,0)$. This complimentary tableau has the area given by the $\sum_{i = 0}^\gamma i |S_{\gamma-i}|$ term. The green squares ($\lambda$) sum to $7$, and the white ones (complement of $\lambda$) sum to $14$.  Putting them together yields the area $\gamma n = 3\cdot 7 = 21$.
\begin{figure}
	\begin{tikzpicture}[every node/.style={minimum size=.5cm-\pgflinewidth, outer sep=0pt}]
    \draw[step=0.5cm,color=black] (-1,-2.5) grid (0.5,1);
    \node[fill=green] at (-0.75,+0.75) {};
    \node[fill=green] at (-0.25,+0.75) {};
    \node[fill=green] at (+0.25,+0.75) {};
    \node[fill=green] at (-0.75,+0.25) {};
    \node[fill=green] at (-0.25,+0.25) {};
    \node[fill=green] at (-0.75,-0.25) {};
    \node[fill=green] at (-0.75,-0.75) {};
    \node[] at (2,0.75) {$0|S_3| = 0\cdot 1 = 0$};
    \node[] at (2,0.25) {$1|S_2| = 1\cdot 1 = 1$};
    \node[] at (2,-0.5) {$2|S_1| = 2\cdot 2 = 4$};
    \node[] at (2,-1.75) {$3|S_0| = 3\cdot 3 = 9$};
\end{tikzpicture}
\caption{Tableau corresponding to $\lambda$ (green) and its complement (white)}\label{fig:tableau}
\end{figure}
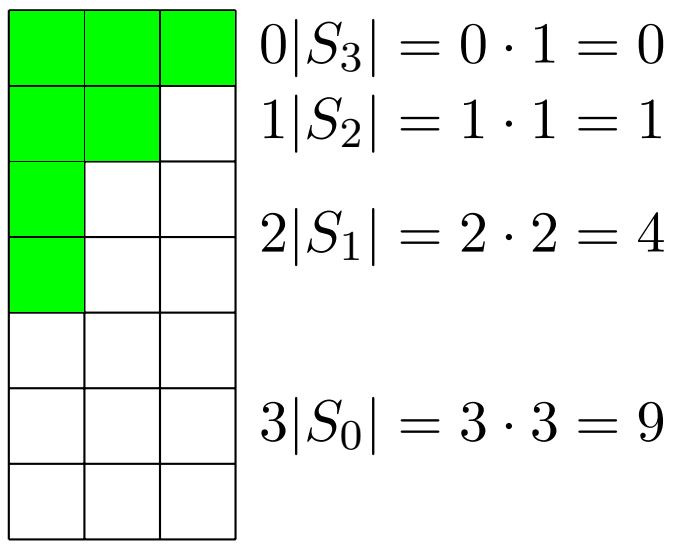
\end{ex}

\begin{thm}\label{thm:correct}
	Let $\cm$ be a linear compartmental model, and let $\lambda$ be such that $\sum_{i=1}^n\lambda_i=n$. If $M_\lambda$ is invertible, then $\bfh_\lambda$ is a parametric profile for $\cj_{\cm}$, and Algorithm \ref{alg} yields the corresponding projections $(f_i,g_k: i \in [n], k \in Out)$, 
	 Where $(g_k :\; k \in Out)$ are the input-output equations, $f_i = x_i-y_i$ for $i \in Out$, and $f_i$ for $i \notin Out$ are given by the pivots in $\bar{M}$ corresponding to $\bfx_{S_0}^{(0)}$.
\end{thm}
\begin{proof}
	First, we need to prove $\bfh_\lambda$ is a parametric profile. This means we need to show $\left\{\bfy^{<\lambda}, \bfu^{(\infty)}\right\}$ is a transcendence basis for $\cc(\bfa)[\bfx^{(\infty)},\bfy^{(\infty)},\bfu^{(\infty)}]/\cj_{\cm}$.
	
	
	Since $M_\lambda$ is invertible, we know that
	\begin{equation*}
		RREF(\bar{M}) = \bordermatrix{
	&\bfx_{S_\gamma}^{(\gamma)}&\bfx_{S^{\leq}_{ \gamma-1}}^{(\gamma-1)} & \bfx_{S^{\leq}_{ \gamma -2}}^{(\gamma-2)} \cdots & \bfx_{S^{\leq}_{ 1}}^{(1)}& \bfx_{S_{0}}^{(0)} & \vrule & \bfx_{S^{>}_{ \gamma -1}}^{(\gamma-1)}& \cdots & \bfx_{S^{>}_{0}}^{(0)} &\vrule & \bfu^{(\leq\gamma-1)} \cr
	&  & & & & &\vrule & & & &\vrule & \cr
	&& & I_{(\gamma-1)n}& & & \vrule & & \star & &\vrule  &\star\star\cr
	&& & & & &\vrule & &  & &\vrule &\cr
	}
	\end{equation*}
The existence of this reduced row echelon form implies that we can write any polynomial in $\cc(\bfa)[\bfx^{(\infty)},\bfy^{(\infty)},\bfu^{(\infty)}]/\cj_{\cm}$ as a polynomial in variables $\left\{\bfx_{S^{>}_{ \gamma -1}}^{(\gamma-1)}, \cdots  ,\bfx_{S^{>}_{0}}^{(0)}, \bfu^{(\infty)}\right\} = \left\{\bfx^{<\lambda}, \bfu^{(\infty)}\right\}$. To prove that these monomials are algebraically independent, we proceed by contradiction:

 Assume there exists a monic irreducible polynomial $p \in \cj_{\cm}\cap \cc(\bfa)[\bfx^{(<\lambda)},\bfu^{(\infty)}]$. Modulo $\cj_{\cm}$, we can rewrite all derivatives of $\bfx$ in terms of linear forms with support $\bfx$. Thus, the set $\left\{ \bfx, \bfu^{(\infty)}\right\}$ is algebraically dependent modulo $\cj_{\cm}$, which is a contradiction. Finally, note that modulo $\cj_{\cm}$, $\bfx^{(<\lambda)} = \bfy^{(<\lambda)}$. This concludes the proof that $\left\{\bfy^{(<\lambda)}, \bfu^{(\infty)}\right\}$ is algebraically independent modulo $\cj_{\cm}$, and therefore, a transcendence basis for $\cc(\bfa)[\bfx^{(\infty)},\bfy^{(\infty)},\bfu^{(\infty)}]/\cj_{\cm}$.

Second, for our polynomials $f_i$, $g_j$ to be projections, we need them to satisfy the following:
\begin{equation*}
\begin{cases}
	f_i \mbox{ is the generator of } \cj_{\cm} \cap \cc(\bfa)[x_i,\bfy^{(<\lambda)},\bfu^{(\infty)}] & 1\leq i \leq n\\
	g_i \mbox{ is the generator of } \cj_{\cm} \cap \cc(\bfa)[\bfy^{(<\lambda + e_j)},\bfu^{(\infty)}] & j\in Out.\\
\end{cases}
\end{equation*}
	Since these are principal ideals, and each $f_i$, $g_j$ is irreducible, this result follows from construction, as we would only need to check each $f_i$ and $g_i$ are contained in their respective ideals. This follows from the fact that $RREF(\bar{M})$  is in reduced row echelon form.
\end{proof}
The idea of the algebraic independence part of the proof of Theorem \ref{thm:correct} is more direct in the case where we do not have inputs, as it makes the quotient ring finite dimensional. We know that the transcendence degree of $\cc(\bfa)[\bfx^{(\infty)},\bfy^{(\infty)}]/\cj_{\cm}$ is $n$. The RREF of $\bar{M}$ shows we can write any polynomial in the quotient ring in terms of $\left\{\bfy^{(<\lambda)}\right\}$. Since $\lambda$ is a weak composition of $n$, we have that $\left\{\bfy^{(<\lambda)}\right\}$ has $n$ elements, thus, must form a transcendence basis for $\cc(\bfa)[\bfx^{(\infty)},\bfy^{(\infty)}]/\cj_{\cm}$.
\begin{ex}\label{ex:ex2} Let us consider the linear compartmental model $\mathcal{M} =  (G, In, Out, Leak)$ with $In = \emptyset$, $Out = \{2,3\}$ and $Leak = [3]$, and $G$ as illustrated in Figure \ref{fig:ex2} (same as Example \ref{ex:ex1} but without the input).
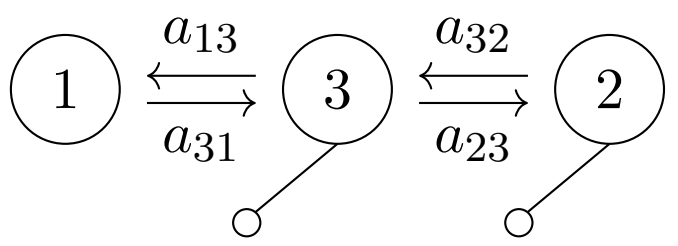
\begin{figure}
\begin{center}
	\begin{tikzpicture}[scale=1.8]
 	\draw (0,0) circle (0.2);	
 	\draw (1,0) circle (0.2);	
 	\draw (2,0) circle (0.2);	
    	\node[] at (0, 0) {1};
    	\node[] at (1, 0) {$3$};
    	\node[] at (2, 0) {$2$};
	 \draw[<-] (0.3, 0.05) -- (0.7, 0.05);	
	 \draw[->] (0.3, -0.05) -- (0.7, -0.05);	
	 \draw[<-] (1.3, 0.05) -- (1.7, 0.05);	
	 \draw[->] (1.3, -0.05) -- (1.7, -0.05);	
   	 \node[] at (0.5, 0.2) {$a_{13}$};
   	 \node[] at (1.5, 0.2) {$a_{32}$};
   	 \node[] at (0.5, -0.2) {$a_{31}$};
   	 \node[] at (1.5, -0.2) {$a_{23}$};
 	\draw (5/3,-.49) circle (0.05);	
	 \draw[-] (2, -.2 ) -- (1.7, -.45);	
 	\draw (2/3,-.49) circle (0.05);	
	 \draw[-] (1, -.2 ) -- (0.7, -.45);	
	\end{tikzpicture}
\end{center}
\caption{Linear compartmental model in Example \ref{ex:ex2}}\label{fig:ex2}
\end{figure}
	Recall $RREF(\bar{M})$ in the model from Example \ref{ex:ex1}:
	\begin{equation*}
	RREF(\bar{M}) = \bordermatrix{
	&\color{blue}y_3''&\color{red}x_1'&\color{blue}y_2'&\color{red}x_1&y_3'&y_2&y_3\cr
&1 & 0 & 0 & 0 & -a_{11}-a_{33} & a_{11}a_{32}-a_{22}a_{32} & a_{11} a_{33}-a_{13} a_{31}-a_{23} a_{32}\cr
& 0 & 1 & 0 & 0 & -\frac{a_{11}}{a_{31}} & \frac{a_{11} a_{32}}{a_{31}} & \frac{a_{11} a_{33}-a_{13} a_{31}}{a_{31}} \cr
 &0 & 0 & 1 & 0 & 0 & -a_{22} & -a_{23}\cr
 &0 & 0 & 0 & 1 & -\frac{1}{a_{31}} & \frac{a_{32}}{a_{31}} & \frac{a_{33}}{a_{31}} \cr
}
\end{equation*}
In this case, $\left\{ \bfy^{(<(0,1,2)} \right\} = \left\{ y_2,y_3, y_3' \right\}$ is a transcendence basis for $\cc(\bfa)[\bfx^{(\infty)},\bfy^{(\infty)}]/\cj_{\cm}$. The associated parametric profile is $\mathbf{h} = (0,0,0,1,2)$, with the corresponding projections:
\begin{align*}
	f_1 =& x_1 -\frac{1}{a_{31}}y_3' + \frac{a_{32}}{a_{31}} y_2 + \frac{a_{33}}{a_{31}} y_3 &\in \cj_{\cm} \cap \cc(\bfa)[x_1,y_2,y_3, y_3']\\
	f_2 =& x_2 - y_2 &\in \cj_{\cm} \cap \cc(\bfa)[x_2,y_2,y_3, y_3']\\
	f_3 =& x_3 - y_3 &\in \cj_{\cm} \cap \cc(\bfa)[x_3,y_2,y_3, y_3']\\
	g_2 =& y_2'-a_{22} y_{2}-a_{23} y_{3} &\in \cj_{\cm} \cap \cc(\bfa)[y_2,y_3,y_2', y_3']\\
	g_3 =& y_{3}''+ \left(-a_{11}-a_{33}\right) y_{3}'+  \left(a_{11}a_{32}-a_{22}a_{32}\right) y_{2}+ \\& \left(a_{11} a_{33}-a_{13} a_{31}-a_{23} a_{32}\right) y_{3}&\in \cj_{\cm} \cap \cc(\bfa)[y_2,y_3, y_3',y_3'']\\
\end{align*}
Note that the structure of $f_1,f_2,f_3$ guarantee that the model is observable, 
as these polynomials allow us to express $x_1$, $x_2$, and $x_3$ in terms of the output variables $y_2$, $y_3$, and their derivatives. 
This means that, while we do not measure $x_1$, for example, we can still compute its value provided we have measurements of $y_3'$, $y_3$, and $y_2$.
\end{ex}

\begin{ex}\label{ex:lindep}
Let us consider the linear compartmental model $\mathcal{M} =  (G, In, Out, Leak)$ with $In = \{1\}$, $Out = \{3\}$ and $Leak = [3]$, and $G$ as illustrated in Figure \ref{fig:lindep}, with $\lambda = (0,0,3)$. 
	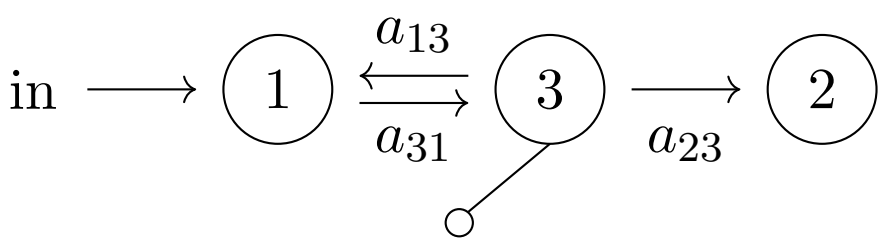
\begin{figure}
	\begin{center}
	\begin{tikzpicture}[scale=1.8]
 	\draw (0,0) circle (0.2);	
 	\draw (1,0) circle (0.2);	
 	\draw (2,0) circle (0.2);	
    	\node[] at (0, 0) {1};
    	\node[] at (1, 0) {$3$};
    	\node[] at (2, 0) {$2$};
	 \draw[<-] (0.3, 0.05) -- (0.7, 0.05);	
	 \draw[->] (0.3, -0.05) -- (0.7, -0.05);	
	 \draw[->] (1.3, 0) -- (1.7, 0);	
   	 \node[] at (0.5, 0.2) {$a_{13}$};
   	 \node[] at (0.5, -0.2) {$a_{31}$};
   	 \node[] at (1.5, -0.2) {$a_{23}$};
 	\draw (2/3,-.49) circle (0.05);	
	 \draw[-] (1, -.2 ) -- (0.7, -.45);	
	 \draw[->] (-0.7, 0) -- (-0.3, 0);	
   	 \node[] at (-.9, 0) {in};
	\end{tikzpicture}
\end{center}
\caption{Linear compartmental model in Example \ref{ex:lindep}}\label{fig:lindep}
\end{figure}
The row reduced echelon form of $\bar{M}$ is given by:
\small
\begin{equation*}
\bordermatrix{
	&y_3'''&x_1''&x_2''&x_1'&x_2'&x_1&x_2&y_3''&y_3'&y_3\cr
&1 & 0 & 0 & 0 & 0 & 0 & 0 & 0 & -a_{11}^{2}-a_{11} a_{33}-a_{13} a_{31}-a_{33}^{2} & \left(a_{11}+a_{33}\right) \left(a_{11} a_{33}-a_{13} a_{31}\right) \cr
& 0 & 1 & 0 & 0 & 0 & 0 & 0 & 0 & \frac{-a_{11}^{2}-a_{13} a_{31}}{a_{31}} & \frac{a_{11} \left(a_{11} a_{33}-a_{13} a_{31}\right)}{a_{31}} \cr
& 0 & 0 & 1 & 0 & 0 & 0 & -a_{22}^{2} & 0 & -a_{23} & -a_{22} a_{23} \cr
&0 & 0 & 0 & 1 & 0 & 0 & 0 & 0 & -\frac{a_{11}}{a_{31}} & \frac{a_{11} a_{33}-a_{13} a_{31}}{a_{31}} \cr
& 0 & 0 & 0 & 0 & 1 & 0 & -a_{22} & 0 & 0 & -a_{23} \cr
&0 & 0 & 0 & 0 & 0 & 1 & 0 & 0 & -\frac{1}{a_{31}} & \frac{a_{33}}{a_{31}} \cr
& 0 & 0 & 0 & 0 & 0 & 0 & 0 & 1 & -a_{11}-a_{33} & a_{11} a_{33}-a_{13} a_{31}\cr
}
\end{equation*}
\normalsize
In this case, $M_\lambda$ consists of the first $7$ columns of the matrix above, and hence the determinant of $M_\lambda$ is 0. 
There is a linear dependence between $y_3'',y_3',y_3$ which can be see from the last row of the matrix above, corresponding to the
equation  $y_3'' + (-a_{11}-a_{33})y_3' + (a_{11} a_{33}-a_{13} a_{31})y_3$.
 Hence,  $\{y_3'',y_3',y_3 \}$  cannot form a transcendence basis for $\cc(\bfa)[\bfx^{(\infty)},\bfy^{(\infty)}]/\cj_{\cm}$. 
\end{ex}

In Example \ref{ex:lindep}, the matrix $M_\lambda$  is not invertible, as the determinant of $M_\lambda$ is 0.  Our goal in the next section is to investigate the cases when $M_\lambda$ is not
of full rank.  We are interested in answering the following question: 
given a linear compartmental model, and a target profile $\bfh_\lambda$, what are the conditions on the graph so that $M_\lambda$ is non-invertible for generic choices of parameters?


\section{Invertibility of $M_\lambda$.}\label{sec:invertibility}

The goals of this section are to explore the combinatorial nature of the
expansion of $\det M_\lambda$  and, in particular, identify criterion that
would guarantee that $\det M_\lambda$ is zero or nonzero.
We associated to $M_\lambda$ a digraph $\cg_\lambda$, and show how the monomials
that appear in the expansion of $\det M_\lambda$ are related to systems
of vertex disjoint paths in $\cg_\lambda$ (Theorem \ref{thm:det}).  This is akin to classic results on
expansions of determinants of matrices related to paths in graphs, like the Gessel-Viennot-Lindstr\"om
Lemma  \cite{GesselViennot1985,Lindstrom1973}.    In addition, using Menger's theorem, we provide a sufficient graphical criterion for $M_\lambda$ to be singular (Theorem \ref{thm:applymengers}).   We also conclude with conjectures related to finding a converse to Theorem \ref{thm:applymengers}.

To illustrate some of the main ideas we will develop, we begin by looking at examples
of the expansion of  $\det(M_\lambda)$ into monomials.

\begin{ex}\label{ex:2case}
Let us consider the linear compartmental model $\mathcal{M} =  (G, In, Out, Leak)$ with $G$ given by the complete graph on 3 vertices, $In = \{1\}$, and the following:
 \begin{enumerate}
		\item $Out = \{2,3\}$, and $\lambda = (0,1,2)$.
		In this case, $$\det(M_\lambda) = a_{31}.$$
		This means that for $\bfh_\lambda$ to be a parametric profile that yields input-output equations with orders given by $\lambda$, there needs to exist an edge from the state variable, $x_1$, to the highest order output variable, $x_3$.
		\item $Out = \{3\}$, and $\lambda = (0,0,3)$.
		In this case, $$\det(M_\lambda) = a_{11} a_{31} a_{32}-a_{31}^{2} a_{12}+a_{21} a_{32}^{2}-a_{22} a_{31} a_{32}.$$
		Rewriting terms we get:
		$$\det(M_\lambda) =  a_{\underline{31}} a_{\underline{11}} a_{32}-a_{\underline{31}} a_{\underline{12}} a_{31}+ a_{\underline{32}}a_{\underline{21}}a_{32}-a_{\underline{32}}a_{\underline{22}} a_{31} .$$
		Note that the underlined terms correspond to paths of length 2 from state variables to the output variable. The non-underlined terms correspond to paths of length 1 from state variables to the output variable. Every monomial in $\det(M_\lambda)$ is represented by 
		$$(\mbox{path $i\rightarrow 3$ of length 2})\cdot(\mbox{path $j\rightarrow 3$ of length 1}),$$
		where $i,j \in \{1,2\}$, $i\neq j$. Thus, if there exists a path of length 2 and a path of length 1 from the set $\{1,2\}$ (state) to $\{3\}$ (output), we can guarantee $M_\lambda$ is not identically 0. We will formalize this condition next.
	\end{enumerate}
\end{ex}

\begin{defn} Given a linear compartmental model, $\cm$, and a target parametric profile $\bfh_\lambda$, we define the following subsets of the set $[n]\times [\gamma]$:
	\begin{align*}
	V_i &= [n] \times \{i\}\\
	L_i &= \{(j,i): \lambda_j < i\} \subseteq V_i \\
	C_i &= \{(j,i): \lambda_j = i\}\subseteq V_i\\
	U_i &= \{(j,i): \lambda_j > i\} \subseteq V_i.
	\end{align*}
\end{defn}

The first coordinate of a pair in $[n]\times [\gamma]$ to corresponds 
to a compartment in $\cm$, while the second coordinate is a grading, 
dependent on the order of derivatives imposed by $\lambda$.

\begin{ex}
Consider the linear compartmental model $\mathcal{M} =  (G, In, Out, Leak)$ with $In = \{1\}$, $Out = \{4,5,6\}$, $Leak = [6]$ and consider the case where $\lambda = (0,0,0,1,2,3)$.   We have:
\begin{equation*}
	\begin{array}{ccc}
		L_1 = \{(1,1),(2,1),(3,1)\} & L_2 = \{(1,2),(2,2),(3,2), (4,2)\} & L_3 = \{(1,3),(2,3),(3,3),(4,3),(5,3)\}\\
		C_1 = \{(4,1)\} & C_2 = \{(5,2)\} & C_3 = \{ (6,3)\}\\
		U_1 = \{(5,1),(6,1)\} & U_2 = \{((6,2))\} & U_3 = \emptyset
	\end{array}.
\end{equation*}
The set $L_1$ is in correspondence with state variables, and we will construct a graph on vertices $[n]\times [\gamma]$, where we can construct paths of prescribed lengths as in Example \ref{ex:2case} that correspond to terms in the determinant of $M_\lambda$.
\end{ex}

\begin{defn}
	Given a linear compartmental model $\cm$ and a target parametric profile $\bfh_\lambda$, 
	we associate a \emph{multipartite auxiliary graph $\cg_\lambda = (V,E)$} 
	with vertices 
	$$V = \bigcup_{i \in [\gamma] }V_i.$$
	The edge set $E$ is defined as follows:
	$$E = \bigcup_{i\in [\gamma-1]} E_{i+1 \leftarrow i}$$
	where
	$$E_{i \leftarrow j} = \{((k,j),(l,i)): k\leftarrow l \in E_G, \; (k,j) \in L_j,\; (l,i) \in V_i\}.$$
\end{defn}

At each one of the $V_i$, we draw an outgoing edge from $(l,i)$ to $(j,i+1)$ if and only if $l \in L_i$ and there was an edge from $l$ to $j$ in the original graph $G$. Schematically, we can see $\cg_\lambda$ is constructed as in Figure \ref{fig:block}. 
\begin{center}
\begin{figure}
	\begin{tikzpicture}
\node[rectangle,
    draw,
    minimum width = 1cm, 
    minimum height = 1cm] (L1) at (0,0) {$L_1$};
\node[rectangle,
    draw,
    align = center,
    minimum width = 1cm, 
    minimum height = 2cm] (C1) at (0,-1.5) {$C_1$\\$U_1$};

\node[rectangle,
    draw,
    minimum width = 1cm, 
    minimum height = 1cm] (L2) at (-3,0) {$L_2$};
\node[rectangle,
    draw,
    align = center,
    minimum width = 1cm, 
    minimum height = 2cm] (C2) at (-3,-1.5) {$C_2$\\$U_2$};
\node (E21) at (-1.5,1) {$E_{2,1}$};
\draw     (-0.5,0.5) -- (-2.5,0.5);
\draw     (-0.5,-0.5) -- (-2.5,-2.5);
\draw[->] (-0.7,-0.25) -- (-2.3,-1.8);
\draw[->] (-0.7,0.25) -- (-2.3,0.25);
\draw[->] (-0.7,0) -- (-2.3,-0.8);

\node[rectangle,
    draw,
    minimum width = 1cm, 
    minimum height = 1cm] (L3) at (-6,0) {$L_3$};
\node[rectangle,
    draw,
    align = center,
    minimum width = 1cm, 
    minimum height = 2cm] (C3) at (-6,-1.5) {$C_3$\\$U_3$};
\node (E32) at (-4.5,1) {$E_{3,2}$};
\draw     (-3.5,0.5) -- (-5.5,0.5);
\draw     (-3.5,-0.5) -- (-5.5,-2.5);
\draw[->] (-3.7,-0.25) -- (-5.3,-1.8);
\draw[->] (-3.7,0.25) -- (-5.3,0.25);
\draw[->] (-3.7,0) -- (-5.3,-0.8);

\node[rectangle,
    draw,
    minimum width = 1cm, 
    minimum height = 1cm] (Lg1) at (-9,0) {$L_{\gamma-1}$};
\node[rectangle,
    draw,
    align = center,
    minimum width = 1cm, 
    minimum height = 2cm] (Cg1) at (-9,-1.5) {$C_{\gamma-1}$\\$U_{\gamma-1}$};

\node[rectangle,
    draw,
    minimum width = 1cm, 
    minimum height = 1cm] (Lg) at (-12,0) {$L_\gamma$};
\node[rectangle,
    draw,
    align = center,
    minimum width = 1cm, 
    minimum height = 2cm] (Cg) at (-12,-1.5) {$C_\gamma$};
\node (Egg1) at (-10.5,1) {$E_{\gamma, \gamma-1}$};
\draw     (-9.5,0.5) -- (-11.5,0.5);
\draw     (-9.5,-0.5) -- (-11.5,-2.5);
\draw[->] (-9.7,-0.25) -- (-11.3,-1.8);
\draw[->] (-9.7,0.25) -- (-11.3,0.25);
\draw[->] (-9.7,0) -- (-11.3,-0.8);
\node (dot) at (-7.5,-1) {$\cdots$};
 
\end{tikzpicture}\caption{Block structure of $\cg_\lambda$.}\label{fig:block}
\end{figure}
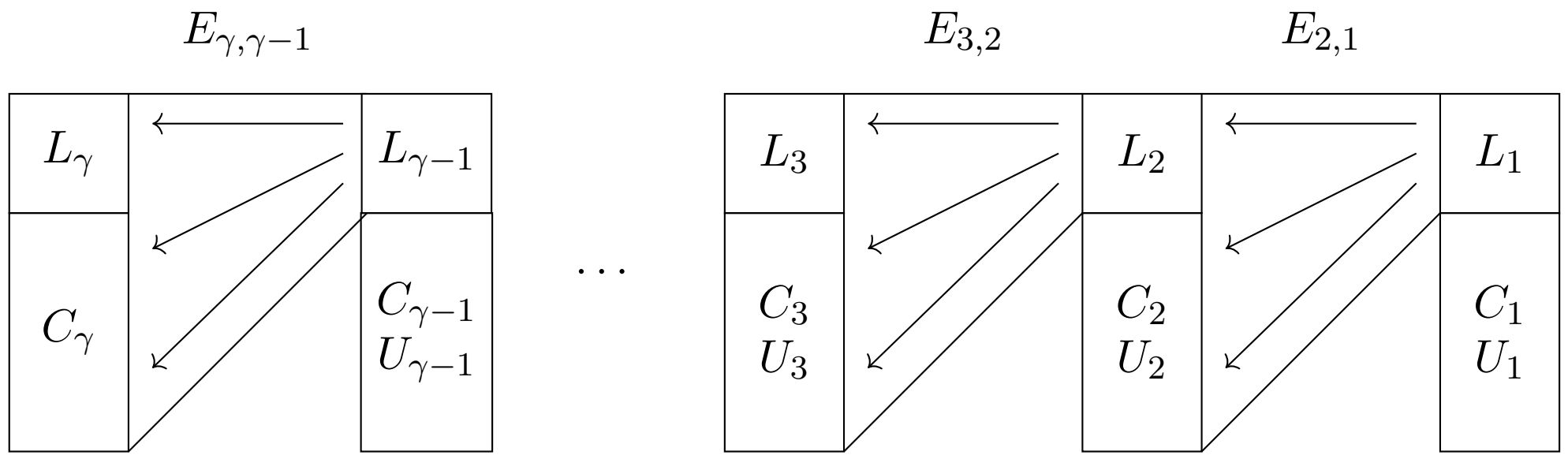
\end{center}

	\begin{figure}
	\begin{center}
		\begin{tikzpicture}[scale=1.8]
			\draw (0,0) circle (0.2);	
			\draw (1,0) circle (0.2);	
			\draw (2,0) circle (0.2);	
			\node[] at (0, 0) {1};
			\node[] at (1, 0) {$2$};
			\node[] at (2, 0) {$3$};
			\draw[->] (0.3, 0) -- (0.7, 0);	
			\draw[->] (1.3, 0) -- (1.7, 0);	
			\node[] at (0.5, 0.2) {$a_{21}$};
			\node[] at (1.5, 0.2) {$a_{32}$};
			\draw (5/3,-.49) circle (0.05);	
			\draw[-] (2, -.2 ) -- (1.7, -.45);	
			\draw[->] (-0.7, 0) -- (-0.3, 0);	
			\node[] at (-.9, 0) {in};
		\end{tikzpicture}
	\end{center}
	\caption{Linear compartmental model in Example \ref{ex:mpgrph}}\label{fig:mpgrph}
\end{figure}
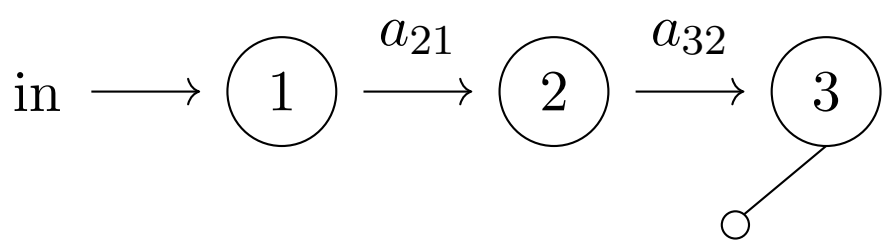

\begin{center}
	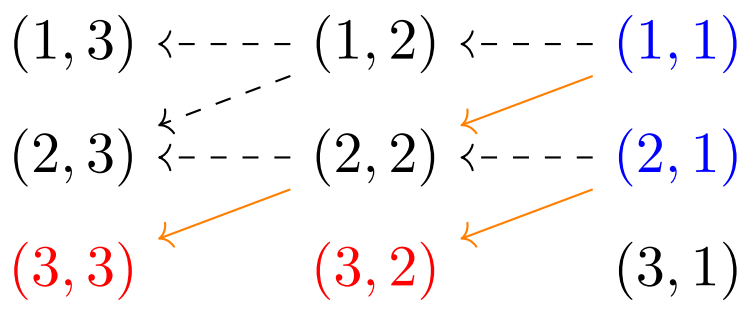
\begin{figure} 
		\begin{tikzpicture}[scale=.5]
			
			\node (11) at (0,1.5) {\color{blue}$(1,1)$};
			\node (21) at (0,0) {\color{blue}$(2,1)$};
			\node (31) at (0,-1.5) {$(3,1)$};
			
			\node (12) at (-4,1.5) {$(1,2)$};
			\node (22) at (-4,0) {$(2,2)$};
			\node (32) at (-4,-1.5) {\color{red}$(3,2)$};
			
			\node (13) at (-8,1.5) {$(1,3)$};
			\node (23) at (-8,0) {$(2,3)$};
			\node (33) at (-8,-1.5) {\color{red}$(3,3)$};
			
			\draw[dashed,->] (11)--(12);
			\draw[->, color = orange] (11)--(22);
			
			\draw[dashed,->] (12)--(13);
			\draw[dashed,->] (12)--(23);
			
			\draw[dashed,->] (21)--(22);
			\draw[->, color = orange] (21)--(32);
			
			\draw[dashed,->] (22)--(23);
			\draw[->, color = orange] (22)--(33);
		\end{tikzpicture}
		\caption{$\cg_\lambda$ from Example \ref{ex:mpgrph}}\label{fig:mpgraph}
	\end{figure}
\end{center}

\begin{ex}\label{ex:mpgrph}
	Let us consider the linear compartmental model $\mathcal{M} =  (G, In, Out, Leak)$ where
	$G$ is the directed path of length $3$,  $In = \{1\}$, $Out = \{3\}$ and $Leak = [3]$, 
	as illustrated in Figure \ref{fig:mpgrph}.  Consider the profile $\lambda= (0,0,3)$.

It follows from Example \ref{ex:2case} that $\det(M_\lambda) = a_{32}a_{21}a_{32}$, by setting
to zero any parameter corresponding to an edge that does not appear in $G$. 
Figure \ref{fig:mpgraph} depicts the augmented graph $\cg_\lambda$. 
Note that the only vertex-disjoint collection of paths (solid edges) 
from state variables (blue) to output variables (red) 
corresponds exactly with the only monomial term in the determinant of $M_\lambda$. This phenomenon will be explained in Theorem \ref{thm:det}.
\end{ex}

Our goal in studying $\cg_\lambda$ is to relate monomials that appear in the expansion of
 $\det(M_\lambda)$ to non-intersecting paths from state variables to output variables. 
 Therefore, we are interested in the number of non intersecting paths from $L_1$ (state variables) to $\bigcup_{i=2}^\gamma (C_i \cup U_i)$ (copies of output variables).

\begin{defn} Let $\cm$ be a linear compartmental model, and let $\lambda$ be such that $\sum_{i \in [n]}\lambda_i = n$.
\begin{enumerate}
	\item For a subgraph $H$ of $\cg_\lambda$, we define the \emph{monomial associated to $H$} as:
	$$ z_H = (-1)^{c} \prod_{(k,j)\leftarrow(l,i) \in E_H} a_{kl},$$
	where $c$ is the number of edge crossings in $H$ not happening on a vertex.
	\item For each $m \in [\gamma-1]$, and each pair of sets of vertices $V_1,V_2 \subseteq V_{\cg_\lambda}$, we define \emph{the set of paths of length $m$ from $V_1$ to $V_2$} as
	 $$\cp_{V_2 \leftarrow V_1}^{m} = \left\{ P\subseteq E_{\cg_\lambda}:\; \mbox{ $P$ is a path of length $m$ from some $v_1 \in V_1$ to some $v_2 \in V_2$} \right\}.$$
	 When either $V_i = \{v_i\}$, we will replace $V_i$ by the vertex $v_i$, to simplify notation.
	\item We define \emph{the set of paths from $V_1$ to $V_2$} as $$\cp_{V_2 \leftarrow V_1} = \bigcup_{m=0}^{\gamma-1} \cp_{V_2 \leftarrow V_1}^{m}  $$
\end{enumerate}
	
\end{defn}

\begin{ex}
	Let $\cm$ be the same linear compartmental model as in Example \ref{ex:mpgrph}. Note that
	\begin{align*}
		\cp_{(2,3)\leftarrow(1,1)}= \left\{ (2,3)\leftarrow (1,2)\leftarrow(1,1), (2,3)\leftarrow (2,2)\leftarrow(1,1) \right\}.
	\end{align*}
	A monomial $z_P$ arising from this set is $$z_{(3,3) \leftarrow (2,2) \leftarrow (1,1)} = a_{32}a_{21}.$$
	If $V_1 = \{(1,1),(2,1)\}$ and $V_2 = \{(3,2),(3,3)\}$, then
	\begin{align*}
		\cp_{V_2\leftarrow V_1} = \{(3,3) \leftarrow (2,2) \leftarrow (1,1), (3,3) \leftarrow (2,2) \leftarrow (2,1), (3,2) \leftarrow (2,1)\}.
	\end{align*}
\end{ex}
\begin{ex}
	In Proposition \ref{prop:Tell}, we make use of several matrices of the form $A_{S^{\leq}_{ \ell},S^{\leq}_{ \ell-1}}$, for some $\ell \in [\gamma-1]$. For the linear compartmental model $\cm = (G,In,Out,Leak)$, where $G$ is the complete graph on 7 vertices (this means $A$ is dense), $In = \emptyset$, $Out = \{4,5,6,7\}$, $Leak = [7]$, and $\lambda = (0,0,0,1,1,2,3)$. Then, $\gamma = \max(\lambda) = 3$ and
	\begin{align*}
		A_{S^{\leq}_{ 1},S^{\leq}_{ 0}} &= A_{[5],[3]} = \begin{bmatrix}
			a_{11} & a_{12} & a_{13}\\
			a_{21} & a_{22} & a_{23}\\
			a_{31} & a_{32} & a_{33}\\
			a_{41} & a_{42} & a_{43}\\
			a_{51} & a_{52} & a_{53}\\
		\end{bmatrix}\\
		A_{S^{\leq}_{ 2},S^{\leq}_{ 1}} &= A_{[6],[5]} = \begin{bmatrix}
			a_{11} & a_{12} & a_{13} & a_{14} & a_{15}\\
			a_{21} & a_{22} & a_{23} & a_{24} & a_{25}\\
			a_{31} & a_{32} & a_{33} & a_{34} & a_{35}\\
			a_{41} & a_{42} & a_{43} & a_{44} & a_{45}\\
			a_{51} & a_{52} & a_{53} & a_{54} & a_{55}\\
			a_{61} & a_{62} & a_{63} & a_{64} & a_{65}\\
		\end{bmatrix}.\\
	\end{align*}
	
\end{ex}

\begin{prop}\label{prop:Tell}
	Let $\lambda$ be such that $\sum_{i \in [n]}\lambda_i = n$. For each $\ell \in [\gamma-1]$, we define 
	\begin{equation*}
		T^{(\ell)} = A_{S^{>}_{\ell},S^{\leq}_{ \ell}}A_{S^{\leq}_{ \ell},S^{\leq}_{ \ell-1}}A_{S^{\leq}_{ \ell-1},S^{\leq}_{ \ell-2}}\cdots A_{S^{\leq}_{ 2},S^{\leq}_{ 1}}A_{S^{\leq}_{ 1},S^{\leq}_{ 0}}
	\end{equation*}
	Then $$T^{(\ell)}_{i-|S^{\leq}_{\ell}|,j} = \sum_{P \in \cp^{\ell-1}_{(i,\ell)\leftarrow (j,1)}} |z_P|$$
	
\begin{proof}
	
	The proof follows from a recursive argument on $m$, where we define 
	\begin{equation*}
	 H^{(m)} = A_{i,S^{\leq}_{ \ell}}A_{S^{\leq}_{ \ell},S^{\leq}_{ \ell-1}}A_{S^{\leq}_{ \ell-1},S^{\leq}_{ \ell-2}}\cdots A_{S^{\leq}_{ \ell-m + 1},S^{\leq}_{ \ell-m}}
	 \end{equation*}
	 We claim that $$ H^{(m)}_j = \sum_{P \in \cp^{m+1}_{(i,\ell) \leftarrow (j,\ell - (m - 1))}} \left|z_P\right|, $$
	 for each $j\in S^{\leq}_{ \ell-m}$. The case of $m=0$ corresponds to
	 \begin{equation*}
	 H^{(0)} = A_{i,S^{\leq}_{ \ell}}
	 \end{equation*}
	 where for each $j\in S_\ell$, $\cp^{1}_{(i,\ell) \leftarrow (j,\ell - 1)} $ is either empty, or it contains a single path of one edge $(i,\ell) \leftarrow (j,\ell-1)$ when $i\leftarrow j \in E_G$. Thus, it follows that $H_j^{(0)} = A_{ij}$ for $j\in S_\ell$. Since we are working on a field where the $a_{ij}$ are indeterminates, we make the distinction between $A_{ij}$, and $a_{ij}$ to be that $A_{ij}$ denotes the $(i,j)$-th entry in the matrix $A$, so it could be 0 if there are no edges between $j$ and $i$. 
	 
	 Now assume our claim holds for all natural numbers $<m$. Then, we have that 
	 $$ H^{(m-1)}_j = \sum_{P \in \cp^{m}_{(i,\ell) \leftarrow (j,\ell - (m-1))}} \left|z_P\right|, $$
	 for each $j\in S^{\leq}_{ \ell-m}$. Observe that 
	 $$H^{(m)} = H^{(m-1)} A_{S^{\leq}_{ \ell-(m-1) },S^{\leq}_{ \ell-m}}.$$
	 This means that for each $r \in S^{\leq}_{ \ell-m}$, we have
	 \begin{align*}
	 	H^{(m)}_r &= H^{(m-1)} A_{S^{\leq}_{ \ell-(m-1)},r}\\
	 	 &= \sum_{j \in S^{\leq}_{ \ell-m}} A_{jr} \sum_{P \in \cp^{m}_{(i,\ell) \leftarrow (j,\ell - m)}} \left|z_P\right|\\
	 	 &= \sum_{j \in S^{\leq}_{ \ell-m}}  \sum_{P \in \cp^{m}_{(i,\ell) \leftarrow (j,\ell - m)}}  \left|A_{jr}z_P\right|
	 \end{align*}
Fix $j$, and note that for each $P \in \cp^{m}_{(i,\ell) \leftarrow (j,\ell - m)}$, we have two cases
\begin{align*}
	 \left|A_{jr}z_P\right| = \begin{cases}
		|a_{jr}z_P| & \mbox{if } (j, \ell-m) \leftarrow (r, \ell-m-1) \in E_{\cg_\lambda}\\
		0 & \mbox{ otherwise.}
	\end{cases}
\end{align*}
Alternatively, we may also write
\begin{align*}
	 \left|A_{jr}z_P\right| = \begin{cases}
		|z_{P'}| & \mbox{if } P' = P\cup \{(j, \ell-m) \leftarrow (r, \ell-m-1)\} \in P^{m+1}_{(i,\ell)\leftarrow (r, \ell-m-1)}\\
		0 & \mbox{ otherwise.}
	\end{cases}
\end{align*}
Thus,
$$H^{(m)}_{r} = \sum_{P \in \cp^{m+1}_{(i,\ell) \leftarrow (r,\ell - m - 1)}} \left|z_P\right|, $$
	for each $r\in S^{\leq}_{ \ell-m-1}$. The proof of the Proposition follows from setting $m=\ell$.\end{proof}

\end{prop}

\begin{prop}\label{prop:det} Let $\cm$ be a linear compartmental model, and $\bfh_\lambda$ a target parametric profile such that $M_\lambda$ is square. Then,
\begin{equation*}
	|\det(M_\lambda)| = \sum_{\sigma \in \mathfrak{S}_{|S_0|}}\text{sign}(\sigma) \sum_{\substack{H_1 \in \cp_{\mathcal{I}_1 \leftarrow \sigma(1)}\\\vdots \\H_{|S_0|} \in \cp_{\mathcal{I}_{|S_0|} \leftarrow \sigma(|S_0|)} }} |z_{H_1}|\cdots |z_{H_{|S_0|}}|,
\end{equation*}
where 
\begin{equation*}
\begin{aligned}
	\mathcal{I} = \bigg[\left(\min\left(S^{\geq}_{2}\right), 2\right),
	\left(\min\left(S^{\geq}_{2}\right)+1,1\right),&\ldots,
	\left(\min\left(S^{\geq}_{2}\right)+ \left|S^{\geq}_{2}\right|,2 \right),\\
	\left(\min\left(S^{\geq}_{3}\right),3)\right),
	\left(\min\left(S^{\geq}_{3}\right)+1,3\right),&\ldots,
	\left(\min\left(S^{\geq}_{3}\right)+ \left|S^{\geq}_{3}\right|,3 \right),\\
	&\vdots\\
	\left(\min\left(S^{\geq}_{\gamma}\right),\gamma\right), 
	\left(\min\left(S^{\geq}_{\gamma}\right) + 1, \gamma\right),
	&\ldots,
	\left(\min\left(S^{\geq}_{\gamma}\right)+ \left|S^{\geq}_{\gamma}\right|,\gamma \right)\bigg].
\end{aligned}
\end{equation*}

	  \begin{proof}
	   Let $\cm$ be a linear compartmental model, and $\bfh_\lambda$ a target parametric profile such that $M_\lambda$ is square. Recall that
	  \begin{equation*}
		M_\lambda = \bordermatrix{
	&\bfx_{S_\gamma}^{(\gamma)}&\bfx_{S^{\leq}_{ \gamma-1}}^{(\gamma-1)} & \bfx_{S^{\leq}_{ \gamma -2}}^{(\gamma-2)}& \cdots & \bfx_{S^{\leq}_{ 1}}^{(1)}& \bfx_{S_{0}}^{(0)}\cr
	& -I_{\left|S_\gamma\right|} & A_{S_{\gamma}, S^{\leq}_{ \gamma-1}} & & && \cr
	&& -I_{\left|S^{\leq}_{ \gamma-1}\right|} & A_{[n],S^{\leq}_{ \gamma -2}}&\ddots & & \cr
	&& && & -I_{\left|S^{\leq}_{ 1}\right|}& A_{[n],S_{0}}\cr
	}
	\end{equation*}
	Consider the following block structure:
	\begin{equation*}
	M_\lambda = \left(\begin{array}{cccccccc}
	 -I_{\left|S_\gamma\right|} & A_{S_{\gamma}, S^{\leq}_{ \gamma-1}} & & && && \\
	& -I_{\left|S^{\leq}_{ \gamma-1}\right|} & A_{S^{\leq}_{ \gamma-1},S^{\leq}_{ \gamma -2}}& & & &&\\
	& &\vrule\overline{\underline{A_{S^{>}_{\gamma -1},S^{\leq}_{ \gamma-2}}}}\vrule & & &&& \\
	& &-I_{\left|S^{\leq}_{ \gamma-2}\right|}& A_{S^{\leq}_{ \gamma -2},S^{\leq}_{\gamma-3}} & &&& \\
	& && \vrule\overline{\underline{A_{S^{>}_{ \gamma -2},S^{\leq}_{\gamma-3}}}}\vrule & &&& \\
	& && -I_{\left|S^{\leq}_{ \gamma -3}\right|} & &&& \\
	& &&&\ddots &&& \\
	& &&&& -I_{\left|S^{\leq}_{ 2}\right|} & A_{S^{\leq}_{2},S^{\leq}_{ 1}}& \\
	&&& && & \vrule\overline{\underline{A_{S^{>}_{2},S^{\leq}_{ 1}}}}\vrule & \\
	& &&&& & -I_{\left|S^{\leq}_{ 1}\right|}& A_{S^{\leq}_{ 1},S_{0}}\\
	& &&&& && \vrule\overline{\underline{A_{S^{>}_{ 1},S_{0}}}}\vrule
	\end{array}\right)
	\end{equation*}
	We perform row reductions to partially eliminate the submatrices of the form $A_{S^{>}_{i},S^{\leq}_{ i-1}}$ (matrices in boxes above). We use the identity matrices to eliminate the block right on top of them. Here is the first iteration of eliminating $ A_{S^{>}_{\gamma -1},S^{\leq}_{ \gamma-2}} $. We want to use the identity matrix $-I_{\left|S^{\leq}_{ \gamma -2}\right|}$ to eliminate the block $A_{S^{>}_{\gamma -1},S^{\leq}_{ \gamma-2}}$. So, we multiply the entire row (block row) that has the identity matrix on the left by $ A_{S^{>}_{\gamma -1},S^{\leq}_{ \gamma-2}} $, and then add that result to the row above. This results in the following matrix:
	
	\begin{equation*}
	M_\lambda = \left(\begin{array}{cccccccc}
	 -I_{\left|S_\gamma\right|} & A_{S_{\gamma}, S^{\leq}_{ \gamma-1}} & & && && \\
	& -I_{\left|S^{\leq}_{ \gamma-1}\right|} & A_{S^{\leq}_{ \gamma-1},S^{\leq}_{ \gamma -2}}& & & &&\\
	& && \vrule\overline{\underline{A_{S^{>}_{\gamma -1},S^{\leq}_{ \gamma-2}} A_{S^{\leq}_{ \gamma -2},S^{\leq}_{\gamma-3}}}}\vrule& &&& \\
	& &-I_{\left|S^{\leq}_{ \gamma-2}\right|}& A_{S^{\leq}_{ \gamma -2},S^{\leq}_{\gamma-3}} & &&& \\
	& && A_{S^{>}_{ \gamma -2},S^{\leq}_{\gamma-3}} & &&& \\
	& && -I_{\left|S^{\leq}_{ \gamma -3}\right|} & &&& \\
	& &&&\ddots &&& \\
	& &&&& -I_{\left|S^{\leq}_{ 2}\right|} & A_{S^{\leq}_{2},S^{\leq}_{ 1}}& \\
	&&& && & A_{S^{>}_{2},S^{\leq}_{ 1}}& \\
	& &&&& & -I_{\left|S^{\leq}_{ 1}\right|}& A_{S^{\leq}_{ 1},S_{0}}\\
	& &&&& && A_{S^{>}_{ 1},S_{0}}
	\end{array}\right)
	\end{equation*}
	Then, we use the $-I_{\left|S^{\leq}_{ \gamma -3}\right|}$ to eliminate $A_{S^{>}_{\gamma -1},S^{\leq}_{ \gamma-2}} A_{S^{\leq}_{ \gamma -2},S^{\leq}_{\gamma-3}}$ and continue in this fashion until we reach the rightmost block. Applying this to all other blocks, we obtain the following row-equivalent matrix:
	\begin{equation*}
	\left(\begin{array}{cccccccc}
	 -I_{\left|S_\gamma\right|} & A_{S_{\gamma}, S^{\leq}_{ \gamma-1}} & & && && \\
	& -I_{\left|S^{\leq}_{ \gamma-1}\right|} & A_{S^{\leq}_{ \gamma-1},S^{\leq}_{ \gamma -2}}& & & &&\\
	& && & &&& T^{(\gamma-1)}\\
	& &-I_{\left|S^{\leq}_{ \gamma-2}\right|}& A_{S^{\leq}_{ \gamma -2},S^{\leq}_{\gamma-3}} & &&& \\
	& && & &&& T^{(\gamma-2)}\\
	& && -I_{\left|S^{\leq}_{ \gamma -3}\right|} & &&& \\
	& &&&\ddots &&& \\
	& &&&& -I_{\left|S^{\leq}_{ 2}\right|} & A_{S^{\leq}_{2},S^{\leq}_{ 1}}& \\
	&&& && & & T^{(2)}\\
	& &&&& & -I_{\left|S^{\leq}_{ 1}\right|}& A_{S^{\leq}_{ 1},S_{0}}\\
	&&& && & & T^{(1)}\\
	\end{array}\right)
	\end{equation*}
	where $T^{(i)} = A_{S^{>}_{i},S^{\leq}_{ i-1}}A_{S^{\leq}_{ i-1},S^{\leq}_{ i-2}}A_{S^{\leq}_{ i-2},S^{\leq}_{ i-3}}\cdots A_{S^{\leq}_{ 2},S^{\leq}_{ 1}}A_{S^{\leq}_{ 1},S^{\leq}_{ 0}}$. Reorganizing the rows, we obtain:
	\begin{equation*}
	\tilde{M_\lambda} = \left(\begin{array}{cccccccc}
	 -I_{\left|S_\gamma\right|} & A_{S_{\gamma}, S^{\leq}_{ \gamma-1}} & & && && \\
	& -I_{\left|S^{\leq}_{ \gamma-1}\right|} & A_{S^{\leq}_{ \gamma-1},S^{\leq}_{ \gamma -2}}& & & &&\\
	& &-I_{\left|S^{\leq}_{ \gamma-2}\right|}& A_{S^{\leq}_{ \gamma -2},S^{\leq}_{\gamma-3}} & &&& \\
	& && -I_{\left|S^{\leq}_{ \gamma -3}\right|} & &&& \\
	& &&&\ddots &&& \\
	& &&&& -I_{\left|S^{\leq}_{ 2}\right|} & A_{S^{\leq}_{2},S^{\leq}_{ 1}}& \\
	& &&&& & -I_{\left|S^{\leq}_{ 1}\right|}& A_{S^{\leq}_{ 1},S_{0}}\\
	& &&&& & & T^{(1)}\\
	&&& && & & T^{(2)}\\
	&&& && & & \vdots\\
	& && & &&& T^{(\gamma-2)}\\
	& && & &&& T^{(\gamma-1)}\\
	\end{array}\right)
	\end{equation*}
	Note that $|\det(M_\lambda)| = |\det(\tilde{M_\lambda})| = |\det(T)|$, where
	$$T = \left(\begin{array}{c}
 T^{(1)}\\
T^{(2)}\\
 \vdots\\
T^{(\gamma-2)}\\
 T^{(\gamma-1)}\\
	\end{array}\right)$$
Furthermore, observe that since $T$ is a square matrix of dimensions $|S_0| \times |S_0|$, we can compute the determinant of $T$ by using permutations as follows:

\begin{equation*}
	|\det(M_\lambda)| = \sum_{\sigma \in \mathcal{S}_{|S_0|}}\mbox{sign}(\sigma)\prod_{i =1 }^{| S_0 |} T_{i,\sigma(i)}
\end{equation*}
We will expand each term to show the relationship of the determinant of $M_\lambda$ with the set of $|L_1|$ paths in $\cg_\lambda$ from $L_1$ to $\bigcup_{i=2}^\gamma (C_i \cup U_i)$.  Proposition \ref{prop:Tell} induces a reindexing of the rows of $T$ that helps interpret this result: to each row $i$ of $T^{(\ell)}$, we index it in $T$ with $\left(i+\left|S^{<}_{\ell}\right|,\ell\right)$. This, in turn means that the rows of $T$ will be indexed by the following vertices of $\cg_\lambda$:
\begin{equation*}
\begin{array}{cc}
	\mathcal{I} = \bigg[\left(\min\left(S^{\geq}_{2}\right), 2\right),
	\left(\min\left(S^{\geq}_{2}\right)+1,2\right),\ldots,
	\left(\min\left(S^{\geq}_{2}\right)+ \left|S^{\geq}_{2}\right|,2 \right), & \in (C_2 \cup U_2)\\
	\left(\min\left(S^{\geq}_{3}\right),3)\right),
	\left(\min\left(S^{\geq}_{3}\right)+1,3\right),\ldots,
	\left(\min\left(S^{\geq}_{3}\right)+ \left|S^{\geq}_{3}\right|,3 \right), & \in (C_3 \cup U_3)\\
	\vdots\\
	\left(\min\left(S^{\geq}_{\gamma}\right),\gamma\right), 
	\ldots,
	\left(\min\left(S^{\geq}_{\gamma}\right)+ \left|S^{\geq}_{\gamma}\right|,\gamma \right)\bigg] & \in (C_{\gamma} \cup U_{\gamma})
\end{array}
\end{equation*}
Proposition \ref{prop:Tell} then shows that each $T_{(i,\ell),j}$ is the sum of all of the monomials corresponding to paths in $\cp_{(i,\ell)\leftarrow (j,1)}^\ell$, for every $(j,1)\in L_1$, and every $(i,\ell) \in \bigcup_{r =2}^{\gamma} C_{r} \cup U_{r}$. Therefore, we have the following expansion:
\begin{equation*}
	|\det(M_\lambda)| = \sum_{\sigma \in \mathfrak{S}_{|S_0|}}\mbox{sign}(\sigma) \sum_{\substack{H_1 \in \cp_{\mathcal{I}_1 \leftarrow (\sigma(1),1)}\\\vdots \\H_{|S_0|} \in \cp_{\mathcal{I}_{|S_0|} \leftarrow (\sigma(|S_0|),1)} }} |z_{H_1}|\cdots |z_{H_{|S_0|}}|
\end{equation*}

  \end{proof}
	  
\end{prop}
\begin{ex}\label{ex:cancel}
	Let $\mathcal{M} = (G, In, Out, Leak)$ be the linear compartmental model where $G$ the complete graph on 6 vertices, $In = \emptyset$, $Out = \{5,6\}$, $Leak = [6]$,  and let us take the target parametric profile $\bfh_\lambda$, where  $\lambda = (0,0,0,0,3,3)$. In this case, we have that 
	\begin{equation*}
		\begin{aligned}
			T &= \left[\begin{array}{c}
				A_{\{5,6\},[4]}\\
				A_{\{5,6\},[4]}A_{[4],[4]}
			\end{array}\right]\\
			&= \left[\begin{array}{cccc}
				a_{51}&a_{52}&a_{53}&a_{54}\\
				a_{61}&a_{62}&a_{63}&a_{64}\\
				\sum_{i = 1}^4a_{5i}a_{i1} & \sum_{i = 1}^4a_{5i}a_{i2} & \sum_{i = 1}^4a_{5i}a_{i3} & \sum_{i = 1}^4a_{5i}a_{i4}\\
				\sum_{i = 1}^4a_{6i}a_{i1} & \sum_{i = 1}^4a_{6i}a_{i2} & \sum_{i = 1}^4a_{6i}a_{i3} & \sum_{i = 1}^4a_{6i}a_{i4}
			\end{array}\right]
		\end{aligned}
	\end{equation*}
	Observe that the columns of $T$ tells us which state variable the path leaves from, and the rows correspond to the output variables the path ends on, with the added grading coming from the length of the paths. So,  we have that
	\begin{equation*}
		\mathcal{I} = [(5,2), (6,2),(5,3),(6,3)].
	\end{equation*}
	Then, Proposition \ref{prop:det} tells us we can compute our determinant as follows:
	\begin{equation*}
		\begin{aligned}
			|\det(M_\lambda)| &= \sum_{\sigma \in \mathfrak{S}_{4}}\mbox{sign}(\sigma) \sum_{\substack{H_1 \in \cp_{(5,2) \leftarrow (\sigma(1),1)}\\H_2 \in \cp_{(6,2) \leftarrow (\sigma(2),1)}\\H_3 \in \cp_{(5,3) \leftarrow (\sigma(3),1)}\\H_4 \in \cp_{(6,3) \leftarrow (\sigma(4),1)}}} |z_{H_1}z_{H_2}z_{H_3}z_{H_4}|\\
			& = \sum_{\sigma \in \mathfrak{S}_{4}}\mbox{sign}(\sigma) a_{5\sigma(1)}a_{6\sigma(2)} \left(\sum_{i = 1}^4a_{5i}a_{i\sigma(3)}\right)\left(\sum_{i = 1}^4a_{5i}a_{i\sigma(4)}\right)
		\end{aligned}
	\end{equation*}
	\end{ex}
	\begin{center}
	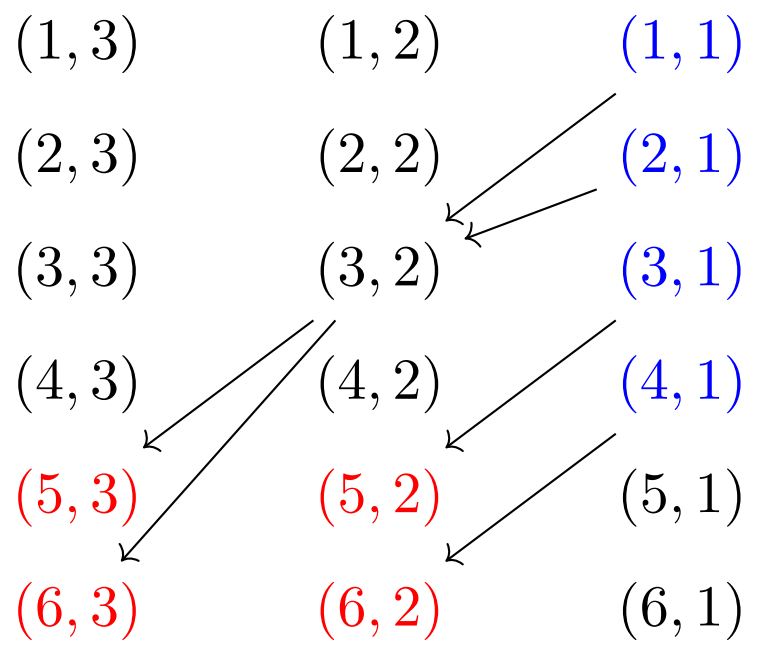
\begin{figure} 
\begin{tikzpicture}[scale=.5]

\node (11) at (0,1.5) {\color{blue}$(1,1)$};
\node (21) at (0,0) {\color{blue}$(2,1)$};
\node (31) at (0,-1.5) {\color{blue}$(3,1)$};
\node (41) at (0,-3) {\color{blue}$(4,1)$};
\node (51) at (0,-4.5) {$(5,1)$};
\node (61) at (0,-6) {$(6,1)$};

\node (12) at (-4,1.5) {$(1,2)$};
\node (22) at (-4,0) {$(2,2)$};
\node (32) at (-4,-1.5) {$(3,2)$};
\node (42) at (-4,-3) {$(4,2)$};
\node (52) at (-4,-4.5) {\color{red}$(5,2)$};
\node (62) at (-4,-6) {\color{red}$(6,2)$};

\node (13) at (-8,1.5) {$(1,3)$};
\node (23) at (-8,0) {$(2,3)$};
\node (33) at (-8,-1.5) {$(3,3)$};
\node (43) at (-8,-3) {$(4,3)$};
\node (53) at (-8,-4.5) {\color{red}$(5,3)$};
\node (63) at (-8,-6) {\color{red}$(6,3)$};

\draw[->] (11)--(32);
\draw[->] (21)--(32);
\draw[->] (32)--(53);
\draw[->] (32)--(63);
\draw[->] (31)--(52);
\draw[->] (41)--(62);
\end{tikzpicture}
\caption{Cancellation occuring in Example \ref{ex:cancel}}\label{fig:cancel}
\end{figure}
\end{center}

	Note that the monomial corresponding to the paths in Figure \ref{fig:cancel} appears twice with opposite sign as it can be obtained from
	\begin{equation*}
		\begin{aligned}
			(5,3) \leftarrow (3,2) \leftarrow \underline{(1,1)}\\
			(6,3) \leftarrow (3,2) \leftarrow \underline{(2,1)}\\
			(5,2) \leftarrow (3,1)\\
			(6,2) \leftarrow (4,1)
		\end{aligned},
	\end{equation*}
	or,
	\begin{equation*}
		\begin{aligned}
			(5,3) \leftarrow (3,2) \leftarrow \underline{(2,1)}\\
			(6,3) \leftarrow (3,2) \leftarrow \underline{(1,1)}\\
			(5,2) \leftarrow (3,1)\\
			(6,2) \leftarrow (4,1)
		\end{aligned}.
	\end{equation*}
	These 2 monomials correspond to the same subgraph (Figure \ref{fig:cancel}), even though they come from different sets of paths. We can see how the terms differ by one transposition, thus they must have opposite signs. This class of cancellations happens because of the fact that the paths that compose the graph are not vertex-disjoint.

	Proposition \ref{prop:det} shows that we can compute the determinant of $M_\lambda$ by looking at paths of certain lengths from the state variables to the outputs. However, there are many cancellations that occur during this process, so we refine our sum even further with the following theorem.
\begin{thm}\label{thm:det}
	Let $\cm = (G, In,Out, Leak)$ be a linear compartmental model, and $\lambda$ such that $M_\lambda$ is square.
	Let $\cn\cp$ be the set consisting of subgraphs of vertex disjoint paths in $\cg_\lambda$ from $L_1$ to $\bigcup_{i=2}^\gamma (C_i \cup U_i)$.
	
	 Then,
	 $$|\det(M_\lambda)|=\sum_{P \in \mathcal{\cn\cp}}z_P.$$
\end{thm}

\begin{proof}
Using Proposition \ref{prop:det}, to prove this theorem, it suffices to show that if there is a set of paths, $P$, that share a vertex, then there exists a set of paths $P'$ that yields the same monomial with opposite sign, and that this process is an involution.  

 Without loss of generality, assume $P$ is composed of paths $H_1,\ldots, H_{|S_0|}$ where $H_i \in \cp_{\mathcal{I}_i \leftarrow (i,1)}$. Furthermore, assume $H_{|S_0|}$, and $H_{\ell}$ share a vertex, say $(r,s)$. Then, let $J_{|S_0|}$ be constructed by taking all the edges from $H_{\ell}$ until we get to $(r,s)$, and then taking all the edges from $H_{|S_0|}$. Let $J_{\ell}$ be constructed by doing the opposite: first taking all the edges from $H_{|S_0|}$ until we get to $(r,s)$, and then taking all the edges from $H_{\ell}$. 
 
 Let $P'$ be composed of $H_1,\ldots, H_{\ell-1},J_\ell, H_{\ell+1}, \ldots, J_{|S_0|}$. Then, note that the monomial corresponding to $P'$ will have different sign as the monomial corresponding to $P$, as it comes from applying the transposition $\tau = (\ell,|S_0|)$ that switches $\ell$ with $|S_0|$. Finally, it follows from construction that the monomials will be equal.
 
 If instead we start with $P'$ and repeat this process, we obtain $K_{|S_0|}$ by taking all the edges from $J_{\ell}$ until we get to $(r,s)$, and then taking all the edges from $J_{|S_0|}$. By construction, $K_{|S_0|}$ is precisely $H_{|S_0|}$. Similarly, $K_{\ell} = H_{\ell}$. Thus, if we consider  $P''$ to be $H_1,\ldots, H_{\ell-1},K_\ell, H_{\ell+1}, \ldots, K_{|S_0|}$, then $P = P''$. Showing this process is an involution.
\end{proof}

While interesting combinatorially, computing the set of all non-intersecting paths of given lengths, and computing determinants using permutations are both extremely computationally expensive. For this reason, we rely on the connection with Menger's theorem is fundamental to make this criterion computationally feasible:
\begin{thm}[Menger's Theorem]
	Let $G= (V,E)$ be a graph and let $S,T \subseteq V$. Then the maximum
number of vertex-disjoint $S-T$ paths is equal to the minimum
size of an $S-T$ separating vertex set.
\end{thm}

The key idea behind using this theorem is that checking the minimum size of an $S-T$ separating vertex set can be done in polynomial time. Thus, we are left with the following result:

\begin{thm}\label{thm:applymengers} 
	Let $\cm = (G, In,Out, Leak)$ be a linear compartmental model, and $\lambda$ such that $M_\lambda$ is square. If the minimum size of a separating set of $L_1$ (identified with $In$) and $\bigcup_{i=2}^\gamma (C_i \cup U_i)$ (identified with $Out$) is less than $|L_1| = |In|$ in $\cg_\lambda$, then the determinant of $M_\lambda$ is zero.
\end{thm}

\begin{proof}
	Assume the minimum size of a separating set of $L_1$ and $\bigcup_{i=2}^\gamma (C_i \cup U_i)$ is less than $|L_1|$. By Menger's Theorem, the maximum number of vertex-disjoint paths between $L_1$ and $\bigcup_{i=2}^\gamma (C_i \cup U_i)$ is strictly less than $|L_1|$. Thus, there is not a single set of vertex-disjoint paths from every element in $L_1$ to $\bigcup_{i=2}^\gamma (C_i \cup U_i)$. By Theorem \ref{thm:det}, $\det(M_\lambda) = 0$.
\end{proof}

We conjecture that the converse of Theorem \ref{thm:applymengers} is true as well. The main obstruction to showing this result directly is that there are still more cancellations happening even after reducing the set we were summing over (Theorem \ref{thm:det}). Here is an example:
\begin{ex}\label{ex:conj}
Let us consider the linear compartmental model $\mathcal{M} =  (G, In, Out, Leak)$ with $In = \{1\}$, $Out = \{4\}$ and $Leak = [4]$, as illustrated in Figure \ref{fig:conjlcm} with $\lambda = (0,0,0,4)$.
	
	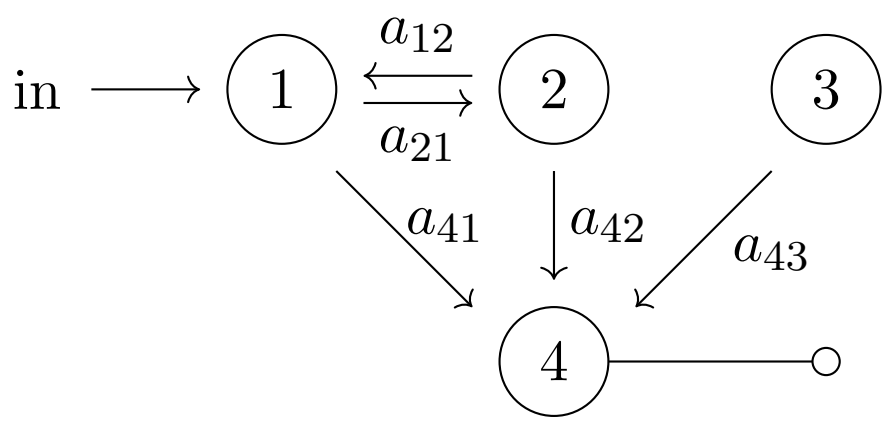
\begin{figure}
	\begin{center}
		\begin{tikzpicture}[scale=1.8]
 	\draw (0,0) circle (0.2);	
 	\draw (1,0) circle (0.2);	
 	\draw (2,0) circle (0.2);	
 	\draw (1,-1) circle (0.2);	
    	\node[] at (0, 0) {1};
    	\node[] at (1, 0) {$2$};
    	\node[] at (2, 0) {$3$};
    	\node[] at (1, -1) {$4$};
	 \draw[<-] (0.3, 0.05) -- (0.7, 0.05);	
	 \draw[->] (0.3, -0.05) -- (0.7, -0.05);	
	 \draw[->] (0.2, -0.3) -- (0.7, -0.8);
	 \draw[->] (1, -0.3) -- (1, -0.7);
	 \draw[->] (1.8, -0.3) -- (1.3, -0.8);
   	 \node[] at (0.5, 0.2) {$a_{12}$};
   	 \node[] at (0.5, -0.2) {$a_{21}$};
   	 \node[] at (0.6, -0.5) {$a_{41}$};
   	 \node[] at (1.2, -0.5) {$a_{42}$};
   	 \node[] at (1.8, -0.6) {$a_{43}$};
 	\draw (2,-1) circle (0.05);	
	 \draw[-] (1.2, -1 ) -- (1.95, -1);	
	 \draw[->] (-0.7, 0) -- (-0.3, 0);	
   	 \node[] at (-.9, 0) {in};
	\end{tikzpicture}
	\end{center}
	\caption{Linear compartmental model in Example \ref{ex:conj}}\label{fig:conjlcm}
	\end{figure}
	Figure \ref{fig:conj} shows two sets of disjoint paths whose monomials cancel each other out in the determinant of $M_\lambda$, and we can see that $G$ is the smallest graph containing these sets of paths. However, we can compute the determinant: 
	$$ \det(M_\lambda) = a_{43}\left(a_{22} a_{11}-a_{11} a_{33}-a_{21} a_{12}-a_{22} a_{33}+a_{33}^{2}\right)\left(a_{11} a_{41} a_{42}-a_{12} a_{41}^{2}+a_{42}^{2} a_{21}-a_{22} a_{41} a_{42}\right) $$
	which shows it is clearly not the zero polynomial, so we conjecture that even if you have a set of paths of the appropriate lengths where cancellation will happen, there must be another non-zero monomial in the determinant. 
	\begin{center}
	\begin{figure} 
\begin{tikzpicture}[scale=.5]

\node (11) at (0,1.5) {\color{blue}$(1,1)$};
\node (21) at (0,0) {\color{blue}$(2,1)$};
\node (31) at (0,-1.5) {\color{blue}$(3,1)$};
\node (41) at (0,-3) {$(4,1)$};

\node (12) at (-4,1.5) {$(1,2)$};
\node (22) at (-4,0) {$(2,2)$};
\node (32) at (-4,-1.5) {$(3,2)$};
\node (42) at (-4,-3) {\color{red}$(4,2)$};

\node (13) at (-8,1.5) {$(1,3)$};
\node (23) at (-8,0) {$(2,3)$};
\node (33) at (-8,-1.5) {$(3,3)$};
\node (43) at (-8,-3) {\color{red}$(4,3)$};

\node (14) at (-12,1.5) {$(1,4)$};
\node (24) at (-12,0) {$(2,4)$};
\node (34) at (-12,-1.5) {$(3,4)$};
\node (44) at (-12,-3) {\color{red}$(4,4)$};

\draw[->,dashed] (11)--(12);
\draw[->,red] (11)--(22);
\draw[->,dashed] (11)--(42);
\draw[->,dashed] (21)--(12);
\draw[->,dashed] (21)--(22);
\draw[->,red] (21)--(42);
\draw[->,red] (31)--(32);
\draw[->,dashed] (31)--(42);

\draw[->,dashed] (12)--(13);
\draw[->,dashed] (12)--(23);
\draw[->,dashed] (12)--(43);
\draw[->,red] (22)--(13);
\draw[->,dashed] (22)--(23);
\draw[->,dashed] (22)--(43);
\draw[->,dashed] (32)--(33);
\draw[->,red] (32)--(43);

\draw[->,dashed] (13)--(14);
\draw[->,dashed] (13)--(24);
\draw[->,red] (13)--(44);
\draw[->,dashed] (23)--(14);
\draw[->,dashed] (23)--(24);
\draw[->,dashed] (23)--(44);
\draw[->,dashed] (33)--(34);
\draw[->,dashed] (33)--(44);

\end{tikzpicture} \qquad
\begin{tikzpicture}[scale=.5]

\node (11) at (0,1.5) {\color{blue}$(1,1)$};
\node (21) at (0,0) {\color{blue}$(2,1)$};
\node (31) at (0,-1.5) {\color{blue}$(3,1)$};
\node (41) at (0,-3) {$(4,1)$};

\node (12) at (-4,1.5) {$(1,2)$};
\node (22) at (-4,0) {$(2,2)$};
\node (32) at (-4,-1.5) {$(3,2)$};
\node (42) at (-4,-3) {\color{red}$(4,2)$};

\node (13) at (-8,1.5) {$(1,3)$};
\node (23) at (-8,0) {$(2,3)$};
\node (33) at (-8,-1.5) {$(3,3)$};
\node (43) at (-8,-3) {\color{red}$(4,3)$};

\node (14) at (-12,1.5) {$(1,4)$};
\node (24) at (-12,0) {$(2,4)$};
\node (34) at (-12,-1.5) {$(3,4)$};
\node (44) at (-12,-3) {\color{red}$(4,4)$};

\draw[->,dashed] (11)--(12);
\draw[->,dashed] (11)--(22);
\draw[->,red] (11)--(42);
\draw[->,red] (21)--(12);
\draw[->,dashed] (21)--(22);
\draw[->,dashed] (21)--(42);
\draw[->,red] (31)--(32);
\draw[->,dashed] (31)--(42);

\draw[->,dashed] (12)--(13);
\draw[->,red] (12)--(23);
\draw[->,dashed] (12)--(43);
\draw[->,dashed] (22)--(13);
\draw[->,dashed] (22)--(23);
\draw[->,dashed] (22)--(43);
\draw[->,dashed] (32)--(33);
\draw[->,red] (32)--(43);

\draw[->,dashed] (13)--(14);
\draw[->,dashed] (13)--(24);
\draw[->,dashed] (13)--(44);
\draw[->,dashed] (23)--(14);
\draw[->,dashed] (23)--(24);
\draw[->,red] (23)--(44);
\draw[->,dashed] (33)--(34);
\draw[->,dashed] (33)--(44);

\end{tikzpicture}
\caption{Cancellation occuring in Example \ref{ex:conj}}\label{fig:conj}
\end{figure}
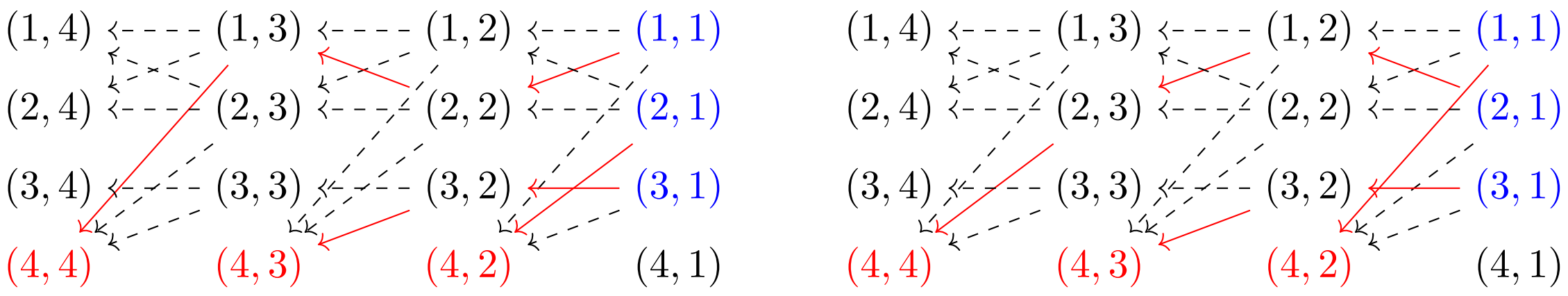
\end{center}
\end{ex}

\begin{ex}\label{ex:conj1}
Consider the linear compartmental model $\mathcal{M} =  (G, In, Out, Leak)$ with $In = \{1\}$, $Out = \{5\}$ and $Leak = [5]$, as illustrated in Figure \ref{fig:conj1lcm} with  $\lambda = (0,0,0,0,5)$.
	\begin{figure}
	\begin{center}
		\begin{tikzpicture}[scale=1.8]
 	\draw (0,0) circle (0.2);	
 	\draw (1,1) circle (0.2);	
 	\draw (2,0) circle (0.2);	
 	\draw (0.5,-1) circle (0.2);	
 	\draw (1.5,-1) circle (0.2);	
    	\node[] at (0, 0) {1};
    	\node[] at (1, 1) {$2$};
    	\node[] at (2, 0) {$3$};
    	\node[] at (0.5, -1) {$5$};
    	\node[] at (1.5, -1) {$4$};
	 \draw[<-] (0.15, 0.25) -- (0.75, 0.85);	
	 \draw[->] (0.25, 0.15) -- (0.85, 0.75);	
	 \draw[->] (0.1, -0.25) -- (0.35, -0.75);
	 \draw[->] (1, 0.7) -- (0.5, -0.7);
	 \draw[->] (1.7, 0) -- (0.3, 0);
	 \draw[->] (1.2, -1) -- (0.8, -1);
 	\draw (-0.5,-1) circle (0.05);	
	 \draw[-] (-0.45,  -1) -- (0.3, -1);	
	 \draw[->] (-0.7, 0) -- (-0.3, 0);	
   	 \node[] at (-.9, 0) {in};
	\end{tikzpicture}
	\end{center}
	\caption{Linear compartmental model in Example \ref{ex:conj1}}\label{fig:conj1lcm}
	\end{figure}
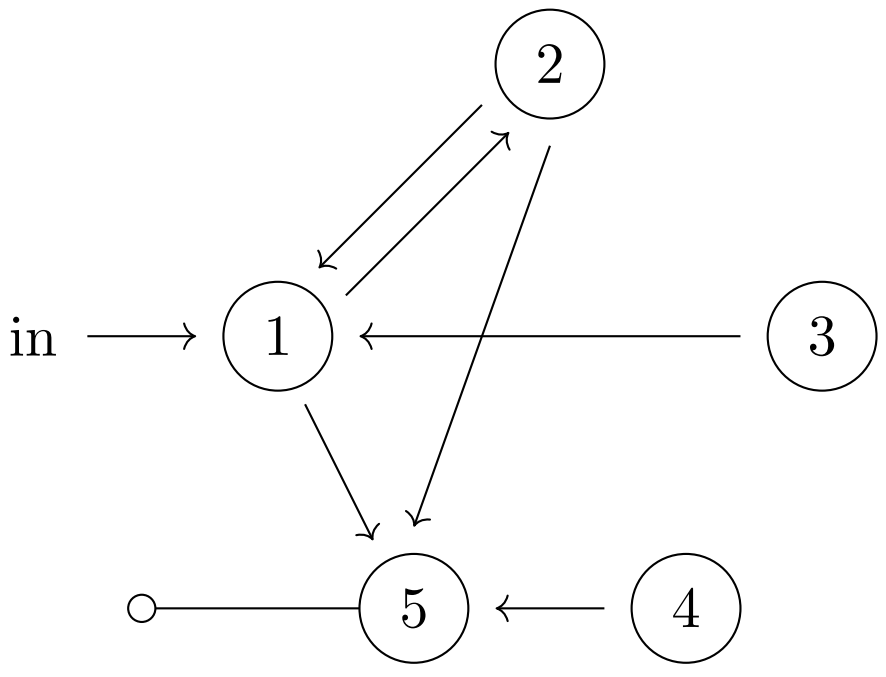
	Figure \ref{fig:conj1} shows two sets of disjoint paths whose monomials cancel each other out in the determinant of $M_\lambda$, and $G$ is the smallest graph containing these sets of paths. However, the determinant of $M_\lambda$ is nonzero. $M_\lambda$ contains 81 distinct monomials.
	\begin{center}
	\begin{figure} 
\begin{tikzpicture}[scale=.5]

\node (11) at (0,1.5) {\color{blue}$(1,1)$};
\node (21) at (0,0) {\color{blue}$(2,1)$};
\node (31) at (0,-1.5) {\color{blue}$(3,1)$};
\node (41) at (0,-3) {\color{blue}$(4,1)$};
\node (51) at (0,-4.5) {$(5,1)$};

\node (12) at (-4,1.5) {$(1,2)$};
\node (22) at (-4,0) {$(2,2)$};
\node (32) at (-4,-1.5) {$(3,2)$};
\node (42) at (-4,-3) {$(4,2)$};
\node (52) at (-4,-4.5) {\color{red}$(5,2)$};

\node (13) at (-8,1.5) {$(1,3)$};
\node (23) at (-8,0) {$(2,3)$};
\node (33) at (-8,-1.5) {$(3,3)$};
\node (43) at (-8,-3) {$(4,3)$};
\node (53) at (-8,-4.5) {\color{red}$(5,3)$};

\node (14) at (-12,1.5) {$(1,4)$};
\node (24) at (-12,0) {$(2,4)$};
\node (34) at (-12,-1.5) {$(3,4)$};
\node (44) at (-12,-3) {$(4,4)$};
\node (54) at (-12,-4.5) {\color{red}$(5,4)$};

\node (15) at (-16,1.5) {$(1,5)$};
\node (25) at (-16,0) {$(2,5)$};
\node (35) at (-16,-1.5) {$(3,5)$};
\node (45) at (-16,-3) {$(4,5)$};
\node (55) at (-16,-4.5) {\color{red}$(5,5)$};

\draw[->,dashed] (11)--(12);
\draw[->,red] (11)--(22);
\draw[->,dashed] (11)--(52);
\draw[->,red] (21)--(12);
\draw[->,dashed] (21)--(22);
\draw[->,dashed] (21)--(52);
\draw[->,dashed] (31)--(12);
\draw[->,red] (31)--(32);
\draw[->,dashed] (41)--(42);
\draw[->,red] (41)--(52);

\draw[->,dashed] (12)--(13);
\draw[->,dashed] (12)--(23);
\draw[->,red] (12)--(53);
\draw[->,dashed] (22)--(13);
\draw[->,red] (22)--(23);
\draw[->,dashed] (22)--(53);
\draw[->,red] (32)--(13);
\draw[->,dashed] (32)--(33);
\draw[->,dashed] (42)--(43);
\draw[->,dashed] (42)--(53);

\draw[->,red] (13)--(14);
\draw[->,dashed] (13)--(24);
\draw[->,dashed] (13)--(54);
\draw[->,dashed] (23)--(14);
\draw[->,dashed] (23)--(24);
\draw[->,red] (23)--(54);
\draw[->,dashed] (33)--(14);
\draw[->,dashed] (33)--(34);
\draw[->,dashed] (43)--(44);
\draw[->,dashed] (43)--(54);

\draw[->,dashed] (14)--(15);
\draw[->,dashed] (14)--(25);
\draw[->,red] (14)--(55);
\draw[->,dashed] (24)--(15);
\draw[->,dashed] (24)--(25);
\draw[->,dashed] (24)--(55);
\draw[->,dashed] (34)--(15);
\draw[->,dashed] (34)--(35);
\draw[->,dashed] (44)--(45);
\draw[->,dashed] (44)--(55);

\end{tikzpicture} \qquad
\begin{tikzpicture}[scale=.5]

\node (11) at (0,1.5) {\color{blue}$(1,1)$};
\node (21) at (0,0) {\color{blue}$(2,1)$};
\node (31) at (0,-1.5) {\color{blue}$(3,1)$};
\node (41) at (0,-3) {\color{blue}$(4,1)$};
\node (51) at (0,-4.5) {$(5,1)$};

\node (12) at (-4,1.5) {$(1,2)$};
\node (22) at (-4,0) {$(2,2)$};
\node (32) at (-4,-1.5) {$(3,2)$};
\node (42) at (-4,-3) {$(4,2)$};
\node (52) at (-4,-4.5) {\color{red}$(5,2)$};

\node (13) at (-8,1.5) {$(1,3)$};
\node (23) at (-8,0) {$(2,3)$};
\node (33) at (-8,-1.5) {$(3,3)$};
\node (43) at (-8,-3) {$(4,3)$};
\node (53) at (-8,-4.5) {\color{red}$(5,3)$};

\node (14) at (-12,1.5) {$(1,4)$};
\node (24) at (-12,0) {$(2,4)$};
\node (34) at (-12,-1.5) {$(3,4)$};
\node (44) at (-12,-3) {$(4,4)$};
\node (54) at (-12,-4.5) {\color{red}$(5,4)$};

\node (15) at (-16,1.5) {$(1,5)$};
\node (25) at (-16,0) {$(2,5)$};
\node (35) at (-16,-1.5) {$(3,5)$};
\node (45) at (-16,-3) {$(4,5)$};
\node (55) at (-16,-4.5) {\color{red}$(5,5)$};

\draw[->,red] (11)--(12);
\draw[->,dashed] (11)--(22);
\draw[->,dashed] (11)--(52);
\draw[->,dashed] (21)--(12);
\draw[->,red] (21)--(22);
\draw[->,dashed] (21)--(52);
\draw[->,dashed] (31)--(12);
\draw[->,red] (31)--(32);
\draw[->,dashed] (41)--(42);
\draw[->,red] (41)--(52);

\draw[->,dashed] (12)--(13);
\draw[->,red] (12)--(23);
\draw[->,dashed] (12)--(53);
\draw[->,dashed] (22)--(13);
\draw[->,dashed] (22)--(23);
\draw[->,red] (22)--(53);
\draw[->,red] (32)--(13);
\draw[->,dashed] (32)--(33);
\draw[->,dashed] (42)--(43);
\draw[->,dashed] (42)--(53);

\draw[->,dashed] (13)--(14);
\draw[->,dashed] (13)--(24);
\draw[->,red] (13)--(54);
\draw[->,red] (23)--(14);
\draw[->,dashed] (23)--(24);
\draw[->,dashed] (23)--(54);
\draw[->,dashed] (33)--(14);
\draw[->,dashed] (33)--(34);
\draw[->,dashed] (43)--(44);
\draw[->,dashed] (43)--(54);

\draw[->,dashed] (14)--(15);
\draw[->,dashed] (14)--(25);
\draw[->,red] (14)--(55);
\draw[->,dashed] (24)--(15);
\draw[->,dashed] (24)--(25);
\draw[->,dashed] (24)--(55);
\draw[->,dashed] (34)--(15);
\draw[->,dashed] (34)--(35);
\draw[->,dashed] (44)--(45);
\draw[->,dashed] (44)--(55);

\end{tikzpicture}
\caption{Cancellation occuring in Example \ref{ex:conj1}}\label{fig:conj1}
\end{figure}
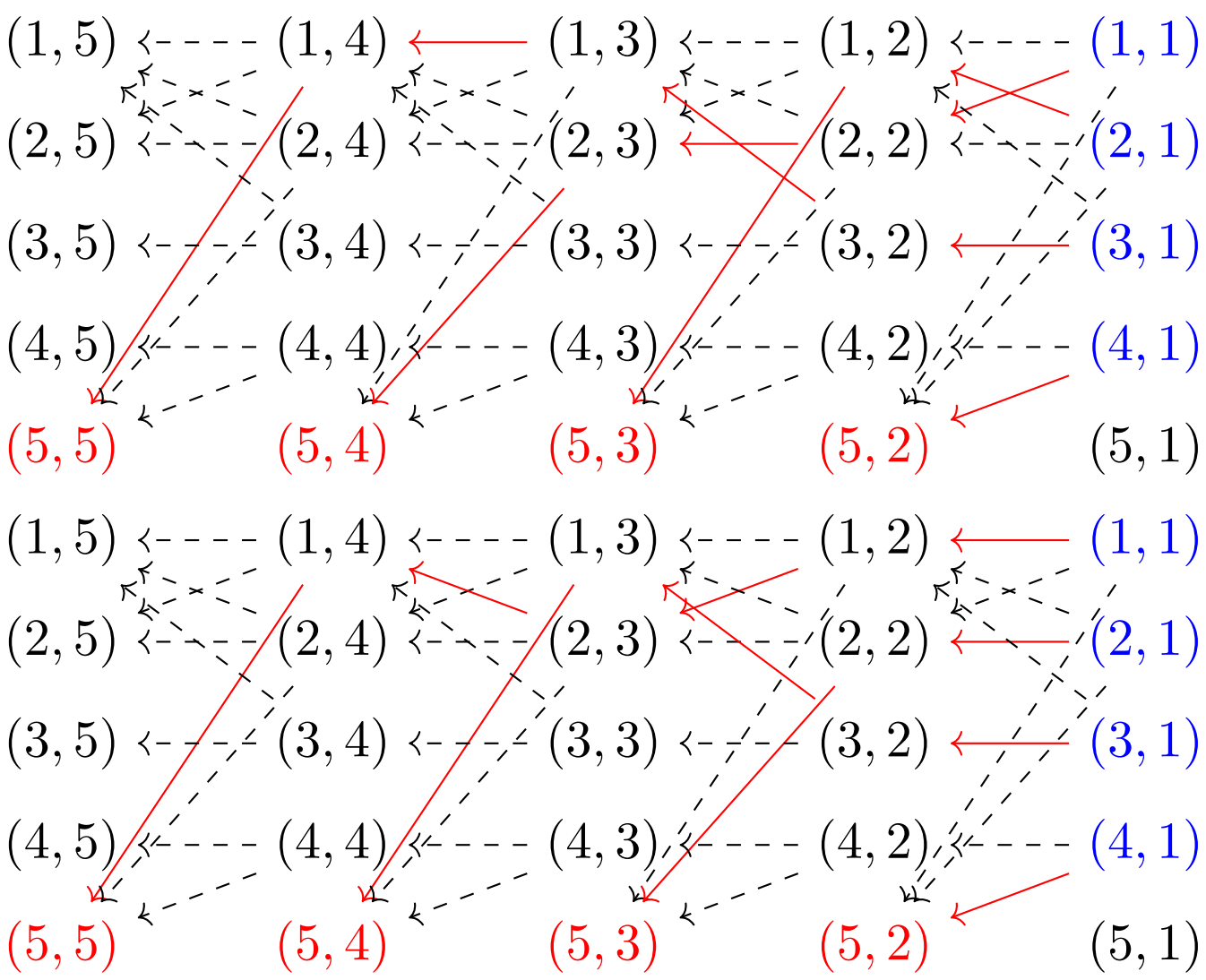
\end{center}
\end{ex}

Examples \ref{ex:conj}  and \ref{ex:conj1} show that  even if there is cancellation of some monomials happening in the formula of $\det(M_\lambda)$ given in Theorem \ref{thm:det}, the existence of a set of vertex-disjoint paths of the appropriate lengths in $\cg_\lambda$ guarantees that $\det(M_\lambda)$ is non-zero. We believe this is due to the fact that we are assuming leaks on every compartment, which in turn gives us edges of the form $(j,i+1)\leftarrow (j,i)$ for every $(j,i)$ in $L_i$, which cannot be deleted from $\cg_\lambda$. The converse of Theorem \ref{thm:det} would follow as a direct consequence of showing that the existence of a set of vertex-disjoint paths of the appropriate lengths is a sufficient condition for $M_\lambda$ to be generically full rank. This would transform our result into an algorithm to check whether input-output equations of a particular order exist.

\section*{Acknowledgments}
Seth Sullivant was partially supported by a Travel Support for Mathematicians Grant 
from the Simons Foundation.

\section*{AI Disclosure Statement}
Artificial intelligence tools were not used in the development of this work.

\bibliography{diffref.bib}{}
\bibliographystyle{plain}

\end{document}